\documentclass{article}

\usepackage{hyperref}
\usepackage[left=2.5cm,right=2.5cm,top=2.5cm,bottom=2.5cm]{geometry}
\usepackage{graphicx}
\usepackage{amssymb,amsmath,amscd}
\usepackage{mathtools, color}
\usepackage{cleveref}
\usepackage{multirow}
\usepackage{soul}

\newcommand{\C}{\mathbb{C}}

\newcommand{\N}{\mathbb{N}}

\newcommand{\R}{\mathbb{R}}
\newcommand{\Z}{\mathbb{Z}}
\newcommand{\Matlab}{\textsc{Matlab}\textsuperscript{\tiny\textregistered}}
\newcommand{\Fref}[1]{\Cref{#1}}%{Figure~\ref{#1}}
\newcommand{\fref}[1]{\cref{#1}}%{Fig.~\ref{#1}}

\newcommand{\eref}[1]{\cref{#1}}
\newcommand{\Sref}[1]{Section~\ref{#1}}

\newcommand{\lp}{\left}
\newcommand{\rp}{\right}

\newtheorem{lem}{Lemma}

\crefformat{equation}{(#2#1#3)}

\title{Pulse-adding of Temporal Dissipative Solitons:\\ Resonant Homoclinic Points and the Orbit Flip of Case B with Delay}

\author{Andrus Giraldo\footnotemark[1] %\footnotemark[3]
\and Stefan Ruschel\footnotemark[2]}% \footnotemark[3]}

\begin{document}

\renewcommand{\thefootnote}{\fnsymbol{footnote}}
\footnotetext[1]{School of Computational Sciences, Korea Institute for Advanced Study, Seoul 02455, Korea

  (\href{mailto:agiraldo@kias.re.kr}{agiraldo@kias.re.kr})}
\footnotetext[2]{Department of Mathematics, The University of
  Auckland, Private Bag 92019, Auckland 1142, New Zealand

  (\href{mailto:stefan.ruschel@auckland.ac.nz}{stefan.ruschel@auckland.ac.nz})}
%\footnotetext[3]{Dodd-Walls Centre for Photonic and Quantum Technologies, New Zealand}
\renewcommand{\thefootnote}{\arabic{footnote}}
\maketitle

%%%%%%%%%%%%%%%%%%%%%%%%%%%%%%%%%%%%%%%%%%%%%%%%%%%%%%%%
\begin{abstract}
  We numerically investigate the branching of temporally localized, two-pulse periodic traveling waves from one-pulse periodic traveling waves with non-oscillating tails in delay differential equations (DDEs) with large delay. Solutions of this type are commonly referred to as temporal dissipative solitons (TDSs) \cite{Yanchuk2019} in applications, and we adopt this term here. We show by means of a prototypical example that \textemdash analogous to traveling pulses in reaction-diffusion partial differential equations (PDEs) \cite{Yanagida1987}\textemdash the branching of two-pulse TDSs from one-pulse TDSs with non-oscillating tails is organized by codimension-two homoclinic bifurcation points  of a real saddle equilibrium \cite{Homburg2010} in a corresponding traveling wave frame. 
  We consider a generalization of Sandstede's model \cite{san1} (a prototypical model for studying codimension-two homoclinic bifurcation points in ODEs) with an additional time-shift parameter, and use \textsc{Auto07p} \cite{Doe1,Doe2} and DDE-BIFTOOL \cite{Sieber2014} to compute numerically the unfolding of these bifurcation points in the resulting DDE. We then interpret this model as the traveling wave equation for TDSs in a DDE with large delay by exploiting the reappearance of periodic solutions in DDEs \cite{Yanchuk2009}. In doing so, we identify both the non-orientable resonant homoclinic bifurcation and the orbit flip bifurcation of case $\mathbf{B}$ as organizing centers for the existence of two-pulse TDSs in the DDE with large delay. Additionally, we discuss how folds of homoclinic bifurcations in an auxiliary system bound the existence region of TDSs in the DDE with large delay.
\end{abstract}
%\keywords{resonant homoclinic orbits, orbit flip t}%Use showkeys class option if keyword
\maketitle

%\tableofcontents
%%%%%%%%%%%%%%%%%%%%%%%%%%%%%%%%%%%%%%%%%%%%%%%%%%%%%%%%

\section{Introduction}\label{sec:Introduction}

The formation of spatiotemporal patterns in nature has captivated researchers for decades; see \cite{Champneys2021} and references therein for a modern review on the topic and methods. Special attention has been paid to localized spatiotemporal phenomena such as traveling waves and pulses, which have great importance, for example, for signal transmission in optical \cite{Lugiato1987} and physiological systems \cite{Hodgkin1952}. Classically, the formation of spatiotemporal patterns is modeled by continuous spatially-extended systems \cite{Cross1993}; on the other hand, spatially-discrete, coarse-grained descriptions in terms of coupled systems have gained more and more attention in the past years \cite{Chow1996}. In both cases, traveling waves and pulses manifest as special solutions with fixed spatial profiles and constant propagation speed corresponding to connecting orbits \cite{Homburg2010} in a suitable co-moving coordinate frame, the so-called traveling wave equation.  
% Crystal growth: Langer1980
% Biological Patterns: GiererA1972Atob
% Instabilities: Cross1993 maybe together with large delay theory. 
% Neural field models: \cite{Coombes}

Recently, such a traveling wave equation has been introduced for systems with time delay \cite{Yanchuk2019}. 
Such time-delay systems, including delay differential equations (DDEs) of the form \begin{equation}\label{eq:DDE}
	u^\prime (t)=f(u(t),u(t-\tau)),
\end{equation}
where $f$ is a sufficiently smooth function in both arguments, are an important class of models for complex dynamical phenomena and the formation of spatiotemporal patterns \cite{Yanchuk2017} in various fields of science \cite{Erneux2009, Smith2010}, including optics \cite{Vladimirov2005} and neuroscience \cite{Laing2022}. Although not explicitly depending  on spatial variables, time-delay systems often arise from feedback loops considering the finite speed of propagation or transport of a physical quantity along a distance. The resulting delayed feedback, when combined with a mechanism of local activation \textemdash for example an excitable configuration \cite{Dubbeldam1999}\textemdash can lead to a regenerative pulsing behavior characterized by a sustained re-excitation repeating approximately every delay interval; see for example \cite{Ruschel2020}.

In analogy to (and to differentiate from) periodic traveling pulses in spatially-extended systems, which are commonly referred to as (dissipative) solitons in physics \cite{Purwins2010, Remoissenet2013}, periodic traveling pulses in DDEs are referred to as \emph{temporal dissipative solitons} (TDSs) \cite{Yanchuk2019}. TDSs have been observed experimentally, for example, semiconductor lasers and optical resonators \cite{garbin2015top, Garbin2020, Marconi2015, romeira2016reg, terrien2020eq}, electrically coupled excitable biological cells \cite{Wedgwood2021}, and various mathematical models of physical processes \cite{Laing2022, Munsberg2020, Pimenov2020, Ruschel2020, Semenov2018, Vladimirov2022, Vladimirov2005}. Of particular interest in spatially extended systems have been the underlying mechanisms that cause the existence of multi-peak traveling pulses \cite{Evans1982, Yanagida1987}. For time-delay systems, multi-peak TDSs have also been observed in models and experiments \cite{ Seidel2022, Terrien2021}; however, their characterization is yet to be as well understood as in the PDE case. To account for this, we will use the traveling wave ansatz for DDEs, developed in \cite{Yanchuk2019}, to identify bifurcations that lead to the branching of two-pulse TDS from one-pulse TDS with non-oscillating tails. 

\subsection{TDSs with non-oscillating tails} 
\begin{figure}
	\centering
	\includegraphics{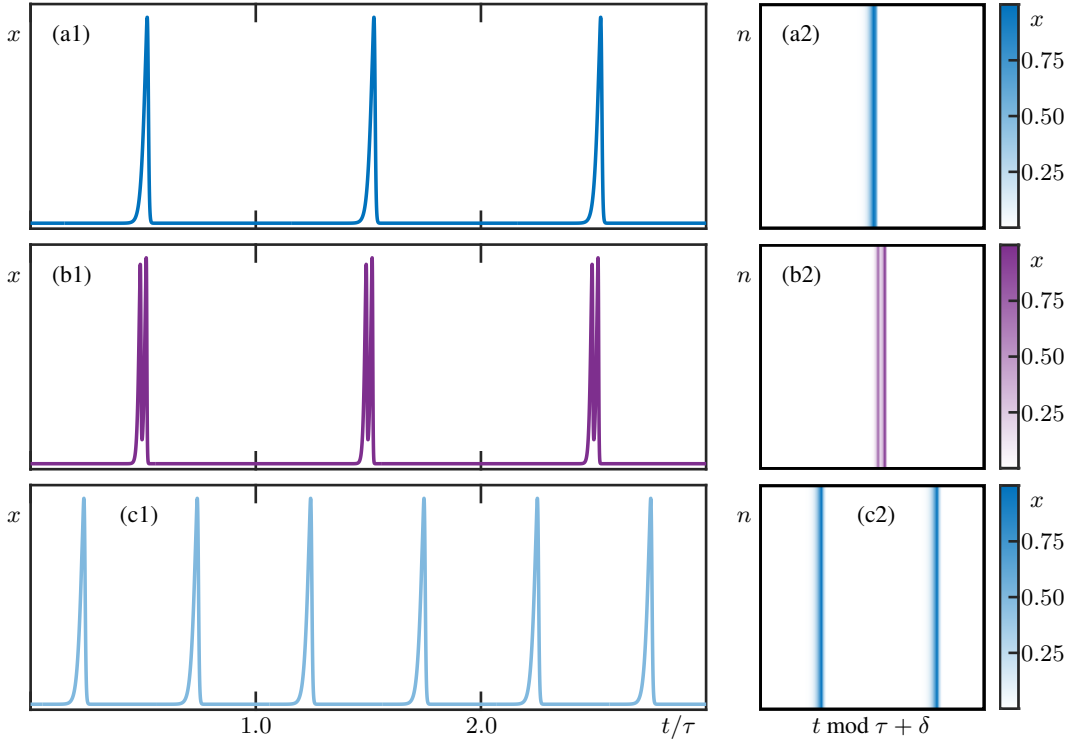}
	\caption{Panels (a1)-(c1) show time traces of temporal dissipative solitons (TDSs) with one, two and two (in order) pulses per delay interval. Panels (a2)--(c2) show space-time plots of the respective TDSs shown in (a1)--(c1). The small correction parameter $\delta$ accounts for the mismatch between the period (two times the period for panels~(c2)) and delay; here, $\delta \approx 0.7033,  0.6605$ and $0.7033$ in panels~(a2), (b2) and (c2), respectively.  Time traces have been obtained from numerical continuation of system \eref{eq:sanDelay} (defined in \Sref{sec:MathModel}). Parameters for all panels are $(a,b,c,\alpha,\gamma,\mu,\tilde{\mu}, \tau) \approx (-1.56,2.5,-1,1,0.5, 0.2,0.03, 100)$. } \label{fig:intSoliton} 
\end{figure}

The focal point of the article is a detailed numerical analysis of a three-component system (denoted $x,y$ and $z$) of DDEs based on a model by Sandstede \cite{san1} and the mathematical framework to view TDSs as approximating homoclinic bifurcations \cite{Lin1986}; see \Sref{sec:MathModel} for the definition and discussion of the model (system \eref{eq:sanDelay}). To illustrate different types of TDSs and the mechanisms that organize them, we present in  \fref{fig:intSoliton} time traces (a1)--(c1) and space-time plots (a2)--(c2) of different types of coexisting periodic solutions for fixed parameters and large delay $\tau$. Panel (a1) shows in blue the $x$-component of a periodic solution with a pronounced, temporally localized peak repeating approximately every delay interval. This localized periodic solution corresponds to a TDS, more specifically, a one-pulse TDS. The analogy to solitons and traveling pulses in spatially-extended systems becomes immediately clear when representing the TDS in a so-called space-time plot. The space-time plot is constructed by splitting the TDS into segments with length equal to period $T$ of the TDS ($2T$ in panel~(c2)) and stacking them vertically, such that the peak appears stationary.  Panel~(a2) shows the corresponding space-time plot of the TDS in panel~(a1) when considering one-hundred periodic segments, i.e, $n=1,2,3,\ldots, 100$. The time along each row (here called pseudo-space) and the ordinal number $n$ (here called pseudo-time) can be seen as analogous to the spatial and temporal variables in spatially extended systems, respectively. This choice of representation is not unique, as there is a certain freedom of choice for the pseudo-spatial variable for delay systems. It is particularly common in optics \cite{Garbin2020, Marconi2015, Terrien2021} to use function segments with delay length $\tau$ resulting in a drift $\delta=T-\tau$ ($\delta=2T-\tau$ for panel~(c2)) of the TDS from each row to the next in the space-time plot. Notice that, by not compensating for the drift $\delta$ in this way, the space-time plot closely resembles the temporal evolution of a traveling pulse with wave speed $\delta$ in a spatially extended system with periodic boundary conditions; thus, motivating the use of the term TDS in the first place. An alternative approach to the spatio-temporal representation of time-delay systems can be found in \cite{Marino2020}, which we do not adopt here. 

%%%%%%%%%%%%%%%%%%%%%%%%%%%%%%%%%%%%%%%%%%%%%%%%%%%%%%%%
%%%%%%%%%%%%%%%%%%%%%%%%%%%%%%%%%%%%%%%%%%%%%%%%%%%%%%%%

Comparing panels (a), (b) and (c), it is an important observation that the one-pulse TDS shown in (a) can coexist with TDSs that have a different pulse shape and, in particular, a different number of pulses per delay interval.  Panels~(b) and (c) show TDSs with two distinct peaks per delay interval, hence referred to as two-pulse TDSs. In order to distinguish between the two-pulse TDS shown in panels (b) and (c), we refer to the TDS in panel (b) as a \emph{bound} two-pulse TDSs (with no exponential localization between the two peaks which are close together) as compared to the \emph{unbound} two-pulse TDS shown  in panel (c) whenever ambiguity arises.

%%%%%%%%%%%%%%%%%%%%%%%%%%%%%%%%%%%%%%%%%%%%%%%%%%%%%%%%
\begin{figure}
	\centering
	\includegraphics{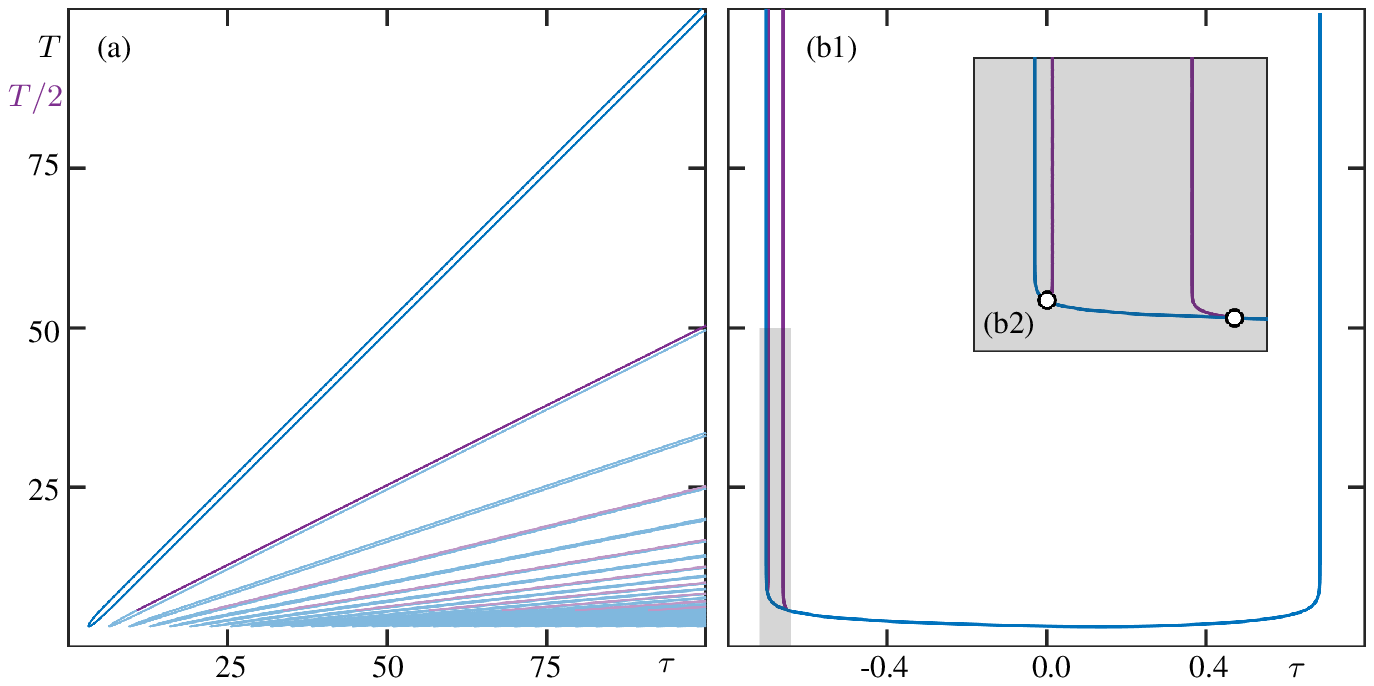}
	\caption{Bifurcation diagram of TDSs from \fref{fig:intSoliton}(a1)-(c2) with respect to the time-delay parameter $\tau$. The curve of one-pulse TDS is shown in vivid blue, whereas the two curves of bound two-pulse TDSs are shown in vivid purple. Unbound TDS with multiple pulses per delay interval are shown in pastel blue, whereas the reappearances of the bound two-pulses curves are shown in pastel purple. Panel (a) shows the period of the orbit (half the period for purple curves) with respect to $\tau$. Panel~(b1) shows the reappearance for $k=-1$ of the vivid curves shown in panel (a). Panel~(b2) shows an enlargement of the grey area in panel~(b1), where white dots indicate period-doubling bifurcations.  Parameters for both panels are $(a,b,c,\alpha,\gamma,\mu,\tilde{\mu}) \approx (-1.56,2.5,-1,1,0.5, 0.2,0.03)$. } \label{fig:intRemergence} 
\end{figure} 
%%%%%%%%%%%%%%%%%%%%%%%%%%%%%%%%%%%%%%%%%%%%%%%%%%%%%%% 

To understand the organization of TDSs in parameter space, we perform a numerical bifurcation analysis.  We show in \cref{fig:intRemergence}(a) the one-parameter $\tau$-bifurcation diagram of the TDS in \fref{fig:intSoliton} for values of $\tau\in[0,100]$. The bifurcation diagram was obtained with the numerical continuation package \textsc{DDE-BIFTOOL} \cite{Sieber2014} for \Matlab. \Fref{fig:intRemergence}(a) shows curves of TDSs (blue, purple, pastel blue, pastel purple) with respect to $\tau$ where the period of the associated TDS is shown on the vertical axis (half of the period for purple curves). Note that the (vivid) blue curve, corresponding to the one-pulse TDS, has two intersection points with the line $\tau=100$. This implies the existence of a secondary one-pulse TDS, which we do not show in \fref{fig:intSoliton}. \Cref{fig:intRemergence}(a) reveals that these two one-pulse TDS are part of the same family of periodic orbits connected through a saddle-node bifurcation for small $\tau$, where the (vivid) blue curve folds back on itself. Notice also that the slope along this family is approximately one, as it is characteristic for one pulse-TDS with scaling $T/\tau \approx 1$ for large $\tau$. 

We remark that \Cref{fig:intRemergence}(a) exhibits a degree of self-similarity that is common to DDEs, where families of periodic orbits with period $T=T(\tau)$ (assuming the family can be locally parameterized by $\tau$) \emph{reappear} under the mapping 
\begin{equation}\label{eq:reap}
\tau\mapsto\tau+kT(\tau)
\end{equation} 
in the $\tau$-bifurcation diagram (shown as pastel blue curves) for any integer $k$ \cite{Yanchuk2009}. In other words, by using the map~\cref{eq:reap} periodic orbits along different families can be identified in the sense that they have exactly the same shape and period despite existing at different values of the delay. For the remainder of the paper, we refer to a family of periodic orbits (and a single periodic orbit) that can be identified under the iteration of map (\ref{eq:reap}) as a \emph{reappearing} family (periodic orbit).
These reappearing families manifest in \cref{fig:intRemergence}(a) as curves that scale with different slope of the form $1/k$ as $\tau$ increases. The phenomenon has several implications, such as a high degree of coexistence (and multi-stability) of periodic orbits for large values of the delay \cite{Yanchuk2017}. 
The effect is clearly visible in \cref{fig:intRemergence}(a), where it becomes virtually impossible to distinguish individual curves with small period for large values of $\tau$.  Now consider a point on the vivid blue curve (corresponding to a one-pulse TDS) in \cref{fig:intRemergence}(a) with sufficiently large $\tau=\tau_0$  with period $T(\tau_0)\approx \tau_0$. Applying the map~\eref{eq:reap} with $k=1$, we obtain a TDS at $\tau_1=\tau_0+T(\tau_0)\approx 2\tau_0$. Comparing the period $T(\tau_1)=T(\tau_0)$ of the so-obtained TDS with $\tau_1$, we have $T(\tau_1)/\tau_1\approx 1/2$, that is, a TDS at $\tau_1$ exhibiting two equidistant peaks per delay interval (twice the period). We have thus discovered a mechanism that generates unbound two-pulse TDSs. More generally,  applying the map~\cref{eq:reap}  generates unbound $(k+1)$-pulse TDSs for any $k \geq 1$. 

We now switch gears and focus our attention on the purple curve in \cref{fig:intRemergence}(a) that corresponds to the bound two-pulse TDSs shown in \cref{fig:intSoliton}(b). Note that the purple curve partially covers the pastel blue curve of unbound two-pulse TDSs. Indeed, the two curves intersect at a single point (left limit of the purple curve) corresponding to a period-doubling bifurcation of the bound two-pulse TDS from the unbound two-pulse TDS. Recall that we chose to plot the purple curve with half of its period such that the period-doubling bifurcation becomes visible as a local bifurcation along the pastel blue branch.

It is worth noting here that the purple branches and the corresponding period-doubling bifurcations reappear for larger values of the delay, but skip every second blue branch. This statement is made precise by the following Lemma; see \Cref{sec:appendix} for the proof. 

\begin{lem}[Reappearance of critical Floquet multipliers]
	\label{thm:reap-bif} 
	Let $\tilde u$ be a periodic solution of Eq.~\eref{eq:DDE} at $\tau=\tilde\tau$ with period $T$. Suppose further that $\tilde u$ possesses a resonant Floquet multiplier $\mu$ such that $\mu^j=1$ for some $j=2,3,\ldots,$ and $\mu^l\neq1$ for all $l=1,2,\ldots,j-1.$	
	
	Then, it holds that $\tilde u$ is a periodic solution of Eq.~\eref{eq:DDE} with a Floquet multiplier $\mu$ for all $\tau=\tilde\tau+jkT$ with $k\in\mathbb{Z}$.
\end{lem}

A direct consequence of \Cref{thm:reap-bif} is the existence of a mapping between points of period-doubling bifurcations (characterized by a critical Floquet multiplier $\mu=-1,$ i.e. $j=2$)
\begin{equation}\label{eq:reap-pd}
	\tau_{\mathbf{PD}}\mapsto\tau_{\mathbf{PD}}+2kT(\tau_{\mathbf{PD}})
\end{equation} 
for every integer $k$. The additional factor two causes the period-doubling to skip every second blue branch in \cref{fig:intRemergence}(a). We have thus established that unbound two-pulse TDSs can be obtained from one-pulse TDSs by using the reappearance map (\ref{eq:reap}), and bound two-pulse TDSs (if present) can be obtained from a period-doubling bifurcation of unbound two-pulse TDSs. We remark that by applying the reappearance map~\eref{eq:reap} to the bound two-pulse TDSs, it is possible to obtain a mixture of unbound and bound TDSs; see pastel purple branches in \fref{fig:intRemergence}(a). The resulting complexity of  reappearing families can be understood by means of the map~\eref{eq:reap} and the connection with underlying homoclinic bifurcations in a suitable traveling wave frame.

\subsection{The traveling wave equation} 
We now focus our attention on the traveling wave equation for TDS, introduced by Yanchuk and colleagues \cite{Yanchuk2019}, that allows us to view TDSs as approximate homoclinic orbits/bifurcations. Notice that the map~\eref{eq:reap} is valid for any integer value of $k$, in particular, for $k=-1$. \Fref{fig:intRemergence}(b) shows the reappearance of the one-pulse TDSs (blue) and the bound two-pulse TDS (purple) for $k=-1$. We remark that the variable $\tau$ along the horizontal axis in panel (b) now corresponds to $-\delta$, where $\delta=T-\tau$ as defined earlier. Therefore, $\tau$ should not be considered a physical parameter in panel~(b), but rather a solution parameter of the system, which can be negative. Similar to the reappearing families in panel (a), each point along the blue and purple curves in \Fref{fig:intRemergence}(b) directly corresponds to TDSs for positive delay. Moreover, \Cref{thm:reap-bif} ensures that certain resonant bifurcations \textemdash in particular period-doubling bifurcations\textemdash     carry over to the bifurcation diagram in panel (b1), where the enlargement (b2)  shows the period-doubling bifurcation points as white dots. 
 
The crucial observation in \Fref{fig:intRemergence} is that in panel (b) the period grows beyond bounds at a finite value of $\tau$, as opposed to panel (a). Such scaling behavior of periodic orbits (here TDS) is typically attributed to a homoclinic bifurcation, where the duration between localized pulses grows without bound. At the point of bifurcation, the TDS coincides with a so-called homoclinic solution, which is a solution of the underlying DDE that is bi-asymptotic to an equilibrium  (also known as the ground state). 
This observation motivates us to call the reparametrization $\bar \tau = \tau - T$ a traveling wave ansatz for TDSs. More specifically, applying the reappearance map~\eref{eq:reap} with $k=-1$ to Eq.~\eref{eq:DDE} and defining $\bar \tau$ as above, we obtain explicitly the DDE (with potentially negative $\bar \tau$)
\begin{equation}\label{eq:DDE-TW}
	u^\prime (t)=f(u(t),u(t-\tau+T))=f(u(t),u(t-\bar \tau)),
\end{equation}
exhibiting homoclinic bifurcations for each branch of TDSs with scaling $T/\tau\approx 1$ as illustrated in \Fref{fig:intRemergence}(b). Hence, we refer to Eq.~\eref{eq:DDE-TW} as the traveling wave equation, or profile equation \cite{Yanchuk2019}.

The traveling wave viewpoint has many advantages and offers additional insight to the analysis of DDEs with large delay. Not only does the bifurcation diagram in panel (b) encode the existence of all reappearing families of TDS in Eq.~\eref{eq:DDE}, but the existence of entire families of TDSs can be conveniently studied by  studying the existence of homoclinic bifurcations in parameter space. 

%%%%%%%%%%%%%%%%%%%%%%%%%%%%%%%%%%%%%%%%%%%%%%%%%%%%%%%%
\begin{figure}
	\centering
	\includegraphics{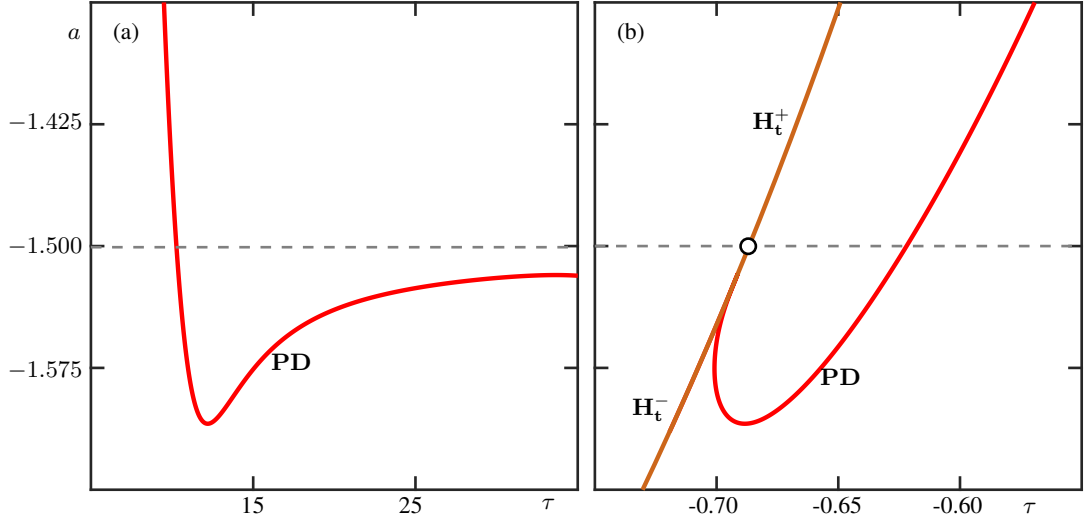}
	\caption{Two-parameter continuation in the $(\tau,a)$-parameter plane of the period-doubling bifurcations in \fref{fig:intRemergence}. Panel (a) shows the continuation (red curve) of the first period-doubling bifurcation $\mathbf{PD}$ in \fref{fig:intRemergence}(a). Panel~(b) shows the continuation (red curve) of the period-doubling bifurcation points $\mathbf{PD}$ \fref{fig:intRemergence}(b2), and homoclinic bifurcation curves $\mathbf{H^-_t}$ and $\mathbf{H^+_t}$ (brown).} \label{fig:intPD} 
\end{figure} 
%%%%%%%%%%%%%%%%%%%%%%%%%%%%%%%%%%%%%%%%%%%%%%%%%%%%%%%

By way of a motivating example, we show in \Fref{fig:intPD}(a)  a two-parameter continuation of the period-doubling bifurcation curve $\mathbf{PD}$ (red) that gives rise to bound two-pulse TDSs in the $(\tau,a)$-parameter plane. The parameter $a$ is the second parameter to be varied; compare \Fref{fig:intRemergence}(a) for the location of the (two) intersection points of the curve $\mathbf{PD}$ with the line $a=-1.56$. \Fref{fig:intPD}(a) suggests that there is an upper bound in $a$ for the curve $\mathbf{PD}$ as $\tau$ increases. The upper limit is interesting as it invariably corresponds to a qualitatively different behavior of the underlying DDE, since for larger values of $a$ there exists only one period-doubling bifurcation in the $\tau$-bifurcation diagram. However, continuing the curve $\mathbf{PD}$ to large values of $\tau$ is challenging as available numerical methods start to break down for large $\tau$ due to the large period and dimension of the underlying boundary value problem.

Let us now change our perspective to the traveling wave viewpoint of TDS. \Fref{fig:intPD}(b) shows the two-parameter bifurcation diagram $(\tau,a)$-parameter plane of the left-most period-doubling bifurcation and homoclinic bifurcation shown in \Fref{fig:intRemergence}(b). Depicted are curves of period-doubling bifurcation $\mathbf{PD}$ (red)  and  homoclinic bifurcations $\mathbf{H_t^+}$ and $\mathbf{H_t^-}$ (brown). Note that the curves $\mathbf{PD}$ in panels (a) and (b) correspond one-to-one by using \cref{thm:reap-bif} and the limiting behavior of the $\mathbf{PD}$ curve in panel (a) can now be clearly identified as being linked to a codimension-two bifurcation along the homoclinic bifurcation curves. Indeed, the codimension-two bifurcation point can be easily identified as resonant homoclinic bifurcation, and thus the upper bound of $\mathbf{PD}$ curve can be identified (dashed black curve). 

\subsection{Codimension-two homoclinic bifurcations of a real saddle}
It is not surprising that homoclinic bifurcations are related to the existence of unbound two-pulse TDSs. Indeed, since the seminal work of Shilnikov \cite{shil5} on the existence of infinitely-many periodic solutions near a homoclinic bifurcation of a saddle focus, countless phenomena relating to the existence/coalescence of periodic solutions, chaotic dynamics, etc. are shown to be organized by homoclinic bifurcations \cite{Homburg2010, Leonid1}. Since then, they have become an invaluable tool for understanding exotic behavior in  applications in biology \cite{Barrio2011, Lina1, Yanagida1987}, electrical and mechanical systems \cite{Champneys1998, Matsumoto1987}, pattern formation \cite{Burke2007, Champneys2021}, optical systems \cite{Bandara2021, Giraldo2021, Parra-Rivas2018, Stitely2022l, Wieczorek2005}, etc. For our purposes, we focus on homoclinic bifurcations to a real saddle, as the TDSs we consider do not have oscillating tails. Unlike their saddle-focus counterpart, homoclinic bifurcations to a real saddle $\mathbf{p}$ only create (as a codimension-one phenomenon) a single periodic solution in three- or higher-dimensional systems of ordinary differential equations. This is true under the following genericity conditions \cite{Homburg2010,Kisaka1993}:
\begin{itemize}
\item[(\textbf{G1})] (Non-resonance) the leading stable $\lambda_s$ and unstable $\lambda_u$ eigenvalues of the real saddle equilibrium are not equal, i.e., $|\lambda_s| \neq \lambda_u$,
\item[(\textbf{G2})] (Principal homoclinic orbit) the homoclinic orbit approaches the equilibrium $\mathbf{p}$ tangent to the weakest stable (unstable) direction of $\mathbf{p}$ in forward (backward) time,
\item[(\textbf{G3})] (Strong inclination) The tangent space of the stable (unstable) manifold of $\mathbf{p}$, followed along the homoclinic orbit in backward (forward) time, does not contain the weak stable (unstable) direction;
\end{itemize}
see \cite{Homburg2010} for a more detailed description of these conditions. By breaking one of the genericity conditions,  the homoclinic bifurcation becomes of codimension-two and more exotic behavior can unfold from it. Particularly, if (\textbf{G1}) is not satisfied, then a resonant bifurcation occurs. The unfolding of the resonant bifurcation depends on orientability properties of the homoclinic orbit, giving rise to an orientable and a non-orientable resonant case \cite{Chow1990}. Particularly, two-homoclinic orbits emanate from the non-orientable resonant bifurcation. On the other hand, if conditions (\textbf{G2}) or (\textbf{G3})  are not fulfilled, then one encounters an orbit flip or inclination flip bifurcation \cite{Deng1993, Hom1, Kisaka1993, Sandstede1993a, Yanagida1987}, respectively. Under additional genericity conditions, both the inclination and the orbit flip bifurcations have three different cases, \textbf{A}, \textbf{B} and \textbf{C},  which are common to both of them, and each case is distinguished by a particular eigenvalue condition of the real saddle \cite{Homburg2010}. The characterizing features of each case are as follows:
\begin{itemize}
\item[\textbf{A.}] A single periodic orbit is created \cite{Kisaka1993, Sandstede1993a}.
\item[\textbf{B.}] The unfolding contains a two-homoclinic bifurcation, a period-doubling bifurcation and a saddle-node bifurcation of periodic orbits \cite{Kisaka1993, Sandstede1993a,Yanagida1987}.
\item[\textbf{C.}] The unfolding contains $k$-homoclinic bifurcations, for any $k \in \N$, and a parameter region where horseshoe and chaotic dynamics exists \cite{Hom1, Nau2002, Sandstede1993a}. 
\end{itemize}

These results concerning the unfoldings of these bifurcations can be extended to DDE case, as it has been shown in \cite{Hupkes2009} for the orbit flip case.

\subsection{Overview of results} We show that, by the reappearance of periodic solutions in DDEs \cite{Yanchuk2009} and the theory of TDSs \cite{Yanchuk2019} developed by Yanchuk and colleagues,  codimension-two homoclinic bifurcations to a real saddle organize the bifurcation diagram of TDSs for large delay. To achieve this, we numerically compute the unfolding of such homoclinic bifurcations and present one-parameter bifurcation diagram near and far-away from these points. In this way, we systematically unveil how different families of TDSs reorganize themselves as different bifurcations  are crossed.  Specifically, we numerically investigate mechanisms that cause bound two-pulse TDSs to branch off from one-pulse TDS. The non-orientable resonant case, and cases~\textbf{B} and \textbf{C} of the orbit and inclination flip bifurcations are natural candidates for our study.  For the purpose of this paper, we restrict ourselves to the codimension-two homoclinic bifurcations that create at most two-pulse TDSs. Thus,  we study the orientable and non-orientable resonant cases, and orbit flip bifurcation of case~\textbf{B}. We leave case~\textbf{C}, which is responsible for higher-order $k$-homoclinic orbits \cite{And2018, Hom1, Sandstede1993a} and  $k$-bound TDSs, for future work. For the codimension-two bifurcations considered in these manuscript, we show that
\vspace{5pt}
\begin{itemize}
	\item[\textbf{(i)}] Unbound two-pulse TDS can be obtained from one-pulse TDS by a mapping argument in parameter space, that is, reappearance of periodic solutions in DDEs \cite{Yanchuk2009}; 
	\item[\textbf{(ii)}] Bound two-pulse TDS can be obtained from unbound two-pulse TDS via period-doubling bifurcations;
	\item[\textbf{(iii)}] The branching of bound two-pulse TDSs from one-pulse TDSs is associated with a period-doubling bifurcation in the traveling wave equation;
	\item[\textbf{(iv)}] The existence of such period-doubling bifurcation is organized by codimension-two homoclinic bifurcation points which give rise to period-doubling bifurcation curves, such as the non-orientable resonant homoclinic bifurcation and the orbit flip bifurcation of type $\mathbf{B}$;
        \item[\textbf{(v)}] TDSs can be created/destroyed by folds of homoclinic bifurcations in an extended traveling wave equation, which are related to saddle-node bifurcations of TDS with infinite period.
\end{itemize}
\vspace{5pt}

We do not delve into the stability properties of TDSs or additional bifurcations emanating for large delay. Here, we focus on their existence as a first step towards their characterization.

\subsection{Organization of the article}
In \Sref{sec:MathModel}, we introduce Sandstede's original model (see system (\ref{eq:san})) and briefly summarize relevant properties pertaining equilibria and homoclinic bifurcations. We showcase the numerical bifurcation diagrams of three codimension-two homoclinic bifurcations that are present, namely, the orientable and non-orientable (twisted) resonant homoclinic bifurcations, and the orbit flip bifurcation of case~\textbf{B}. Following this, we give a detailed description of how we modify Sandstede's original model by introducing a time-shift parameter $\tau$. 
In Secs.~\ref{sec:ResonantCases}--\ref{sec:OrbitFlipCase}, we showcase the numerical unfoldings of the orientable and non-orientable (twisted) resonant homoclinic bifurcations, and the orbit flip bifurcation of case~\textbf{B} in our extended model by lifting the homoclinic bifurcation from system (\ref{eq:san}) to a homoclinic bifurcation at non-zero $\tau$ in the modified system (\ref{eq:sanDelay}) performing a homotopy step in $\tau$.
We present in \Sref{sec:OrientableResonantCase} and \Sref{sec:NonOrientableResonantCase} two-parameter bifurcations diagrams of the orientable and non-orientable resonant homoclinic bifurcations in system (\ref{eq:sanDelay}), respectively, and discuss in detail the resulting one-parameter bifurcation diagrams for TDSs for  large delay.  Particularly for the non-orientable resonant case, we report in \Sref{sec:twoInfty} on the relation between a reduced multivalued-map and the existence of  TDSs.  
In \Sref{sec:OrbitFlipCase}, we present the numerical unfolding of the orbit flip one-parameter bifurcation in system (\ref{eq:sanDelay}) and discuss the resulting large delay bifurcation diagram for TDS.  
We summarize our results and discuss avenues for future research in \Sref{sec:conclusions}.

%%%%%%%%%%%%%%%%%%%%%%%%%%%%%%%%%%%%%%%%%%%%%%%%%%%%%%%%
\section{Mathematical model}\label{sec:MathModel}
Sandstede's model \cite{san1} is a three-dimensional vector field constructed to investigate codimension-two homoclinic bifurcations of a real saddle equilibrium; namely, resonant bifurcations, inclination flip bifurcations and orbit flip bifurcations distinguished by specific parameter choices. Sandstede's model has been employed to understand transitions between different codimension-two homoclinic bifurcations as codimension-three phenomena \cite{Hom2000, OldKra1}, and to characterize the organization of invariant manifolds near them \cite{Agu1, And1, And2018}.

We make use of the existence of these codimension-two homoclinic bifurcations to study their numerical unfoldings by means of a time shift parameter.  In \Sref{sec:NoDelaySan}, we  first review general properties of Sandstede's model and present numerical bifurcation diagrams for the relevant codimension-two homoclinic bifurcation points. In \Sref{sec:intDelaySan}, we explain our choice of time shift parameter and how we lift the numerical results of the ODE case to the DDE case.

\subsection{Resonant and Orbit flip bifurcations in Sandstede's model} \label{sec:NoDelaySan}
Here, we focus on the orientable and non-orientable resonant cases \cite{Chow1990}, and the orbit flip of case~\textbf{B} \cite{And1, Hom1, Homburg2010, Kisaka1993, Sandstede1993a, Yanagida1987}. As such, we consider the following system of equations introduced by Sandstede: 
\begin{equation}\label{eq:san}
\begin{aligned} 
\dot x &= ax + by -ax^2+(\tilde{\mu}-\alpha z)x(2-3x), \\[-0.5mm]
\dot y &= bx +ay -\frac32 bx^2-\frac32 axy-2y(\tilde{\mu}-\alpha
z), \\[-0.5mm]
\dot z &= cz +\mu x +\gamma xz,
\end{aligned}
\end{equation}
where additional parameters $\beta$ and $\delta$ have been set to zero; see \cite{san1} for the full system. We remark that the $z$-axis is invariant under the flow, and the $xy$-plane is also invariant if $\mu=0$. System~\eref{eq:san} was constructed such that, when $(\mu,\tilde{\mu})=(0,0)$ and $0<a^2 \leq b^2$, there exist a homoclinic orbit converging bi-asymptotically to the origin $\mathbf{0}$, which is an equilibrium of system~\eref{eq:san} for any parameter value. This homoclinic orbit is contained in the $xy$-plane and is given by the cartesian leaf
\begin{equation}\label{eq:cartLeaf}
x^2(1-x)-y^2 = 0,
\end{equation}
for positive $x$ and $y$. The associated homoclinic bifurcation persists as a curve in the $(\mu,\tilde{\mu})$-parameter plane; however, the underlying homoclinic orbit is not contained in the $xy$-plane for an arbitrary small perturbation away from $(\mu,\tilde{\mu})=(0,0)$, as shown by Sandstede \cite{san1}. 

The eigenvalues of the origin $\mathbf{0}$ are given by
\begin{equation} \label{eq:EigCond}
\lambda_{1,2}= a \pm \sqrt{b^2+4\tilde{ \mu}^2} \text{ and }
\lambda_3=c,
\end{equation}
where the eigenvector $e_3$ associated with $\lambda_3$ points in the $z$-direction.  For the remainder of this paper, we restrict parameters to $$(b,c,\alpha,\gamma)=(2.5,-1,1,0.5)$$ such that $\lambda_2<\lambda_3<0<\lambda_1$ for $a^2< c^2<b^2+4\tilde{ \mu}^2$, and we denote by $\sigma=\lambda_1+\lambda_3$ the sum of the leading unstable and stable eigenvalues (saddle quantity). In this configuration, $\mathbf{0}$ is a real saddle equilibrium with two-dimensional stable manifold $W^s(\mathbf{0})$ and one-dimensional unstable manifold $W^u(\mathbf{0})$, where $e_3$ corresponds to the weak stable eigenvector. Notice that, when $(\mu,\tilde{\mu})=(0,0)$, the unstable eigenvector $e_1$ associated with $\lambda_1$ and the strong stable direction $e_2$ associated with $\lambda_2$ span the $xy$-plane. Therefore, the homoclinic orbit given by Eq.~\eref{eq:cartLeaf} is in orbit flip configuration \cite{Sandstede1993a, san1} as the homoclinic orbit converges tangentially to the strong stable direction in forward time. In particular, an orbit flip bifurcation with respect to a stable manifold is a codimension-two phenomenon that marks the transition when $W^s(\mathbf{0})$ changes from orientable (homeomorphic to a cylinder locally around the homoclinic orbit) to non-orientable (homeomorphic to a M\"obius band locally around the homoclinic orbit) \cite{Sandstede1993a}. We remark that this is not the only codimension-two homoclinic bifurcation that allows such transition between orientability and non-orientability; namely, inclination flip bifurcations are another mechanism for this transition to occur \cite{Hom1, Homburg2010, Kisaka1993, Sandstede1993a, Yanagida1987}.

In the following, we use \textsc{Auto07p} to compute parameter values where different codimension-two homoclinic bifurcations occur for system~\eref{eq:san}.
%
%%%%%%%%%%%%%%%%%%%%%%%%%%%%%%%%%%%%%%%%%%%%%%%%%%%%%%%%
\begin{figure}
\centering
\includegraphics{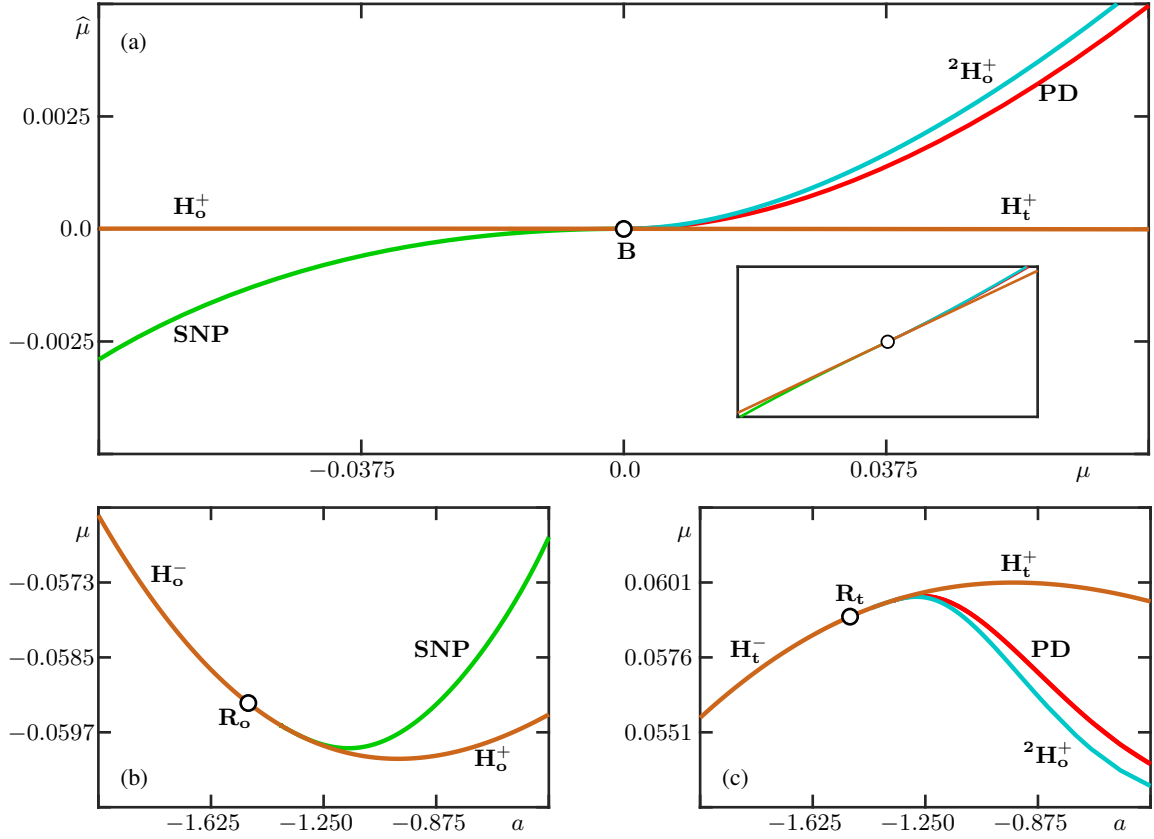}
\caption{Bifurcation diagram in the $(\mu,\widehat{\mu})$-parameter plane of system~\eref{eq:san} near an orbit flip bifurcation point $\mathbf{B}$ of case~$\mathbf{B}$ in panel~(a); here, $\widehat{\mu}=\tilde{\mu} - 0.5047\mu$. The inset of panel~(a) shows the bifurcation diagram in the original $(\mu,\tilde{\mu})$-parameter plane. Panel~(b) shows the bifurcation diagram in the $(a,\mu)$-parameter plane near an orientable resonant bifurcation $\mathbf{R_o}$, and panel~(c) shows the situation near an non-orientable (twisted) resonant homoclinic bifurcation $\mathbf{R_t}$ also in the $(a,\mu)$-parameter plane. The right superscript on the homoclinic bifurcation curves indicates the sign of the saddle quantity.  Shown are curves of homoclinic bifurcations~$\mathbf{H_o}$ and $\mathbf{H_t}$ (brown), saddle-node bifurcations $\mathbf{SNP}$ of periodic orbits (green), period-doubling bifurcations $\mathbf{PD}$ (red), and two-homoclinic bifurcations $\mathbf{H^2_o}$ (cyan). Parameters are $(b,c,\alpha,\gamma)=(2.5,-1,1,0.5)$ for all panels. In panel~(a) $a$ is fixed to $-0.5$, in panel~(b) $\tilde{\mu}$ is fixed to $-3\cdot 10^{-2}$, and in panel~(c) $\tilde{\mu}$ is fixed to  $3\cdot 10^{-2}$. } \label{fig:bifAutoNoDelay} 
\end{figure} 
%%%%%%%%%%%%%%%%%%%%%%%%%%%%%%%%%%%%%%%%%%%%%%%%%%%%%%%%
%
\Fref{fig:bifAutoNoDelay}(a) shows the numerical unfolding, in the $(\mu,\tilde{\mu})$-parameter plane, of an orbit flip bifurcation point $\mathbf{O_B}$ of case~\textbf{B} in system~\eref{eq:san} for $a=-0.5$; see \cite{Kisaka1993, Sandstede1993a, Homburg2010}. For visualization purposes, we choose to present the bifurcation diagram with respect to $\widehat{\mu}=\tilde{\mu} - 0.5047\mu$ on the vertical axis with the inset showing original parameters. Here, the homoclinic bifurcation curve (brown) changes from orientable $\mathbf{H^+_o}$ to non-orientable $\mathbf{H^+_t}$ as $\mu$ changes sign, where the superscript indicates the sign of the saddle quantity $\sigma$. Three bifurcation curves emanate generically from the homoclinic flip bifurcation of case~\textbf{B}; namely, a curve of saddle-node bifurcation of periodic solution $\mathbf{SNP}$ (green), a period-doubling bifurcation curve $\mathbf{PD}$ (red) and an orientable two-homoclinic bifurcation curve $\mathbf{^2H^+_o}$ (cyan); see  \cite{And1} for a detailed characterization of this bifurcation. The curve $\mathbf{^2H^+_o}$  is one of the main objects of study in this manuscript, as each point along $\mathbf{^2H^+_o}$ corresponds to system~\eref{eq:san} exhibiting a homoclinic orbit that pulses two-times before converging biasymptotically to $\mathbf{0}$.

Orbit flip bifurcations are not the only mechanism capable of creating two-homoclinic bifurcation curves. In fact, for three- or higher-dimensional vector fields, so-called resonant bifurcations can also give rise to two-homoclinic bifurcation curve in parameter plane \cite{Homburg2010}. Resonant bifurcations are codimension-two bifurcations with the property that the saddle quantity $\sigma$ changes sign along a codimension-one homoclinic bifurcation curve, that is, $\sigma=0$ at the resonant bifurcation.  The unfolding of a resonant bifurcation depends on wether the homoclinic orbit is orientable or non-orientable \cite{Chow1990, Kisaka1993, Homburg2010}.  To illustrate the unfolding of resonant bifurcations, we choose points $(\mu,\tilde{\mu}) \approx (-0.0595, -0.03)$  and $(\mu,\tilde{\mu})\approx (0.0595, 0.03)$  in  \fref{fig:bifAutoNoDelay}(a), which correspond to points along the curves $\mathbf{H^+_o}$ and  $\mathbf{H^+_t}$, respectively. As a next step, we do continuation of the homoclinic bifurcation curves in the $(a,\mu)$-parameter plane by fixing $\tilde{\mu}$. For the chosen parameter values, that is $\tilde{\mu}=-0.03$ and $\tilde{\mu}=0.03$, $\sigma=0$ implies that  $$a=a_{R} :=-c-\sqrt{b^2+4\tilde{\mu}^2} \approx -1.5007$$ for both the orientable and non-orientable case.  Indeed, we systematically find both resonant cases by making use of the bifurcation diagram shown in \fref{fig:bifAutoNoDelay}~(a). 

\Fref{fig:bifAutoNoDelay}(b) shows the numerical unfolding, in the $(a,\mu)$-parameter plane, of an orientable resonant bifurcation point $\mathbf{R_o}$ in system~\eref{eq:san} for $\tilde{\mu}=-0.03$.  Notice that the saddle quantity along $\mathbf{H^+_o}$ changes at $\mathbf{R_o}$, which is reflected in the change of superscript to $\mathbf{H^-_o}$. From the point  $\mathbf{R_o}$, a saddle-node bifurcation curve of periodic orbits emanates. Notice for the orientable resonant case, no two-homoclinic bifurcation curve is involved in the unfolding. However, this changes for the non-orientable case. Compare  \Fref{fig:bifAutoNoDelay}(c) showing the numerical unfolding of the non-orientable resonant bifurcation point $\mathbf{R_t}$ for $\tilde{\mu}=0.03$, where we find a two-homoclinic bifurcation and period-doubling bifurcation curves emanating from $\mathbf{R_t}$. 

\subsection{Introducing a suitable time shift parameter in Sandstede's model} \label{sec:intDelaySan}
We now introduce a time shift parameter $\tau$ to system~\eref{eq:san} to define the DDE 
\begin{equation}\label{eq:sanDelay} 
\begin{aligned} 
\dot x &= ax + by -ax^2+(\tilde{\mu}-\alpha z)x(2-3x), \\[-0.5mm]
\dot y &= bx +ay -\frac32 bx^2-\frac32 axy-2y(\tilde{\mu}-\alpha
z)+ \kappa \text{S}_{\tau}\! \lp[xy\rp], \\[-0.5mm]
\dot z &= cz +\mu x +\gamma xz,
\end{aligned}
\end{equation}
where $\text{S}_\tau$ is the left-translation operator by $\tau$, i.e. $\text{S}_\tau\!\lp[ \phi \rp](t) = \phi(t-\tau)$ for all $t \in \R$. Notice that the new parameter $\kappa$ controls the influence of the newly introduced time shift parameter, such that we recover system~\eref{eq:san} when $\kappa=0$. By introducing the delay term this way, we take advantage of our numerical results in \cref{sec:NoDelaySan}. More precisely, we use DDE-BIFTOOL to perform a numerical homotopy step in $\kappa$ to lift solutions of the ODE, system~\eref{eq:san}, to the DDE, system~\eref{eq:sanDelay}, for a given $\tau$. 

We specifically choose this delay term to preserve certain properties of system~\eref{eq:san}. Namely, the $z$-axis remains invariant under the flow of system~\eref{eq:sanDelay}, which simplifies to an one-dimensional ODE along this axis. Similarly to system~\eref{eq:san}, the $xy$-plane is invariant when $(\mu,\tilde{\mu})=(0,0)$; however for $\kappa \neq 0$ and $\tau \neq 0$, the dynamics restricted to the $xy$-plane is no longer finite dimensional.  Notice also that the constant solution segment $\mathbf{\widehat{0}}$, defined as $\mathbf{\widehat{0}}(\theta)=(0,0,0)$ for all $\theta \in [-\tau,0]$, is an equilibrium for all parameter values analogous to $\mathbf{0}$ in system~\eref{eq:san}. Since the time shift parameter enters system~\eref{eq:sanDelay} in the nonlinear terms, it does not affect the spectrum of $\mathbf{\widehat{0}}$; that is, its eigenvalues coincide with  the eigenvalue spectrum of $\mathbf{0}$ in system~\eref{eq:san} as given by Eqs.~\eref{eq:EigCond}. 
It follows that homoclinic orbit exhibited by system~\eref{eq:sanDelay} must be in an orbit flip configuration when $(\mu,\tilde{\mu})=(0,0)$. This is due to the eigenvalue configuration and the invariance of the $xy$-plane, as discussed above. This justifies our specific choice for the delay term, and it allows to readily identify orbit flip bifurcations without any additional setup. 

Our primary goal now is to study the numerical unfolding of the orbit flip and the resonant bifurcations in system~\eref{eq:sanDelay} with respect to the time shit parameter $\tau$. A secondary goal is to illustrate and gain additional insight into the orientability of homoclinic orbits in DDEs. As opposed to the ODE case, asserting orientability is not straightforward analytically or even numerically. Thus, we systematically study the unfoldings of these codimension-two points by a priori knowing the orientability properties for the ODE case in system~\eref{eq:san}, which is separated by the orbit flip in parameter space. We start with the resonant bifurcations and move to the orbit flip bifurcation of case~\textbf{B}.

In the following, we find convenient to introduce the extended system
\begin{equation} \label{eq:sanDelayExt}
\begin{aligned} 
\dot x &= ax + by -ax^2+(\tilde{\mu}-\alpha z)x(2-3x), \\[-0.5mm]
\dot y &= bx +ay -\frac32 bx^2-\frac32 axy-2y(\tilde{\mu}-\alpha
z)+ \kappa \text{S}_{\tau}\!\lp[xy\rp], \\[-0.5mm]
\dot z &= cz +\mu x +\gamma xz, \\[-0.5mm]
\dot \tau &= 0,
\end{aligned}
\end{equation}
which corresponds to system~\eref{eq:sanDelay} appended with the trivial equation $\dot \tau = 0$. This notion of extended system will be valuable for characterizing minima and maxima of bifurcation curves in the $(\tau,a)$- and $(\tau, \mu)$-parameter plane of system~\eref{eq:sanDelay}.

\section{Resonant homoclinic bifurcations} \label{sec:ResonantCases}

We begin by studying the orientable and non-orientable cases of resonant homoclinic bifurcations. Since numerical procedures to test orientability for a homoclinic orbit in the DDE case are not readily available, we take advantage of knowing the orientability in the ODEs case. We pick a homoclinic orbit of system~\eref{eq:san} and continue it as a large-period periodic orbit of system~\eref{eq:sanDelay} to moderate values of the delay.

\subsection{Codimension-two orientable resonant case} \label{sec:OrientableResonantCase}
Here, we present a numerical bifurcation diagram of the orientable resonant homoclinic bifurcation in the $(\tau,a)$-parameter plane for fixed $\mu,\tilde{\mu},\kappa$ and \textemdash to remind the reader\textemdash fixed parameters $(b,c,\alpha,\gamma)=(2.5,-1,1,0.5)$ as in \cref{sec:NoDelaySan}. As a starting point, we choose the homoclinic solution of system~\eref{eq:san} at 
$$(a,\mu,\tilde{\mu}) \approx (-0.06, -0.0597, -0.03)$$
along the curve $\mathbf{H^-_o}$ shown in \fref{fig:bifAutoNoDelay}(b). We then perform a series of two-parameter continuation runs to lift the homoclinic orbit from system~\eref{eq:san} to system~\eref{eq:sanDelay} with parameter values 
$$(a,\mu,\tilde{\mu},\kappa,\tau)= (-0.6,-0.2,-0.0550,-1,-1.0558),$$ as summarised in \cref{tab:oriSeqCont}. The thus obtained homoclinic bifurcation serves as a starting point for our numerical bifurcation analysis in the $(\tau,a)$-parameter plane. It is convenient to also lift the saddle-node bifurcation shown in \fref{fig:bifAutoNoDelay}(b) to the same $(\mu,\tilde{\mu},\kappa)$ parameter values; the procedure is largely analogous and is also given in \cref{tab:oriSeqCont}. 
\begin{table}[]
\centering
\begin{tabular}{|ccccccrcl|}
\hline
\multicolumn{5}{|c|}{\textbf{Parameters}} & \multicolumn{1}{c|}{\multirow{2}{*}{\textbf{\begin{tabular}[c]{@{}c@{}}Continuation\\ parameters\end{tabular}}}} & \multicolumn{3}{c|}{\multirow{2}{*}{\textbf{\begin{tabular}[c]{@{}c@{}}Stopping \\ condition\end{tabular}}}} \\ \cline{1-5}
\multicolumn{1}{|c|}{$a$} & \multicolumn{1}{c|}{$\mu$} & \multicolumn{1}{c|}{$\tilde{\mu}$} & \multicolumn{1}{c|}{$\kappa$} & \multicolumn{1}{c|}{$\tau$} & \multicolumn{1}{c|}{} & \multicolumn{3}{c|}{} \\ \hline
\multicolumn{9}{|c|}{\textbf{Homoclinic solution}} \\ \hline
\multicolumn{1}{|c|}{$-0.6$} & \multicolumn{1}{c|}{$-0.0597$} & \multicolumn{1}{c|}{$-0.03$} & \multicolumn{1}{c|}{$0$} & \multicolumn{1}{c|}{$0$} & \multicolumn{1}{c|}{$\kappa,\tilde{\mu}$} & $\kappa$ & $=$ & $-1.0$ \\ \hline
\multicolumn{1}{|c|}{$-0.6$} & \multicolumn{1}{c|}{$-0.0597$} & \multicolumn{1}{c|}{$-0.1987$} & \multicolumn{1}{c|}{$-1.0$} & \multicolumn{1}{c|}{$0$} & \multicolumn{1}{c|}{$\tilde{\mu},\tau$} & $\tilde{\mu}$ & $=$ & $-0.0550$ \\ \hline
\multicolumn{1}{|c|}{$-0.6$} & \multicolumn{1}{c|}{$-0.0597$} & \multicolumn{1}{c|}{$-0.0550$} & \multicolumn{1}{c|}{$-1.0$} & \multicolumn{1}{c|}{$-0.7151$} & \multicolumn{1}{c|}{$\mu,\tau$} & $\mu$ & $=$ & $-0.2$ \\ \hline
\multicolumn{1}{|c|}{$-0.6$} & \multicolumn{1}{c|}{$-0.2$} & \multicolumn{1}{c|}{$-0.0550$} & \multicolumn{1}{c|}{$-1.0$} & \multicolumn{1}{c|}{$-1.0558$} & \multicolumn{1}{c|}{\textemdash} & \multicolumn{3}{c|}{\textemdash} \\ \hline
\multicolumn{9}{|c|}{\textbf{Saddle-node bifurcation of periodic solution}} \\ \hline
\multicolumn{1}{|c|}{$-0.3852$} & \multicolumn{1}{c|}{$-0.0550$} & \multicolumn{1}{c|}{$-0.03$} & \multicolumn{1}{c|}{$0$} & \multicolumn{1}{c|}{$0$} & \multicolumn{1}{c|}{$a,\kappa$} & $a$ & $=$ & $-0.8$ \\ \hline
\multicolumn{1}{|c|}{$-0.8$} & \multicolumn{1}{c|}{$-0.0550$} & \multicolumn{1}{c|}{$-0.03$} & \multicolumn{1}{c|}{$-0.0118$} & \multicolumn{1}{c|}{$0$} & \multicolumn{1}{c|}{$\kappa,\tilde{\mu}$} & $\kappa$ & $=$ & $-1.0$ \\ \hline
\multicolumn{1}{|c|}{$-0.8$} & \multicolumn{1}{c|}{$-0.0550$} & \multicolumn{1}{c|}{$-0.1990$} & \multicolumn{1}{c|}{$-1.0$} & \multicolumn{1}{c|}{$0$} & \multicolumn{1}{c|}{$\tilde{\mu},\tau$} & $\tilde{\mu}$ & $=$ & $-0.0550$ \\ \hline
\multicolumn{1}{|c|}{$-0.8$} & \multicolumn{1}{c|}{$-0.0550$} & \multicolumn{1}{c|}{$-0.0550$} & \multicolumn{1}{c|}{$-1.0$} & \multicolumn{1}{c|}{$-0.7285$} & \multicolumn{1}{c|}{$\mu,\tau$} & $\mu$ & $=$ & $-0.2$ \\ \hline
\multicolumn{1}{|c|}{$-0.8$} & \multicolumn{1}{c|}{$-0.2$} & \multicolumn{1}{c|}{$-0.0550$} & \multicolumn{1}{c|}{$-1.0$} & \multicolumn{1}{c|}{$-1.1037$} & \multicolumn{1}{c|}{\textemdash} & \multicolumn{3}{c|}{\textemdash} \\ \hline
\end{tabular}\vspace{6pt}
\caption{Additional information for the computation of the bifurcation diagram for the orientable resonant case in system~\eref{eq:sanDelay}. Shown are the sequences of continuation runs to lift the homoclinic bifurcation and bifurcations of periodic solutions to non-zero $\tau$. Other parameters are given in \Sref{sec:NoDelaySan}.}\label{tab:oriSeqCont}
\end{table}

%%%%%%%%%%%%%%%%%%%%%%%%%%%%%%%%%%%%%%%%%%%%%%%%%%%%%%%%
\begin{figure}
\centering
\includegraphics{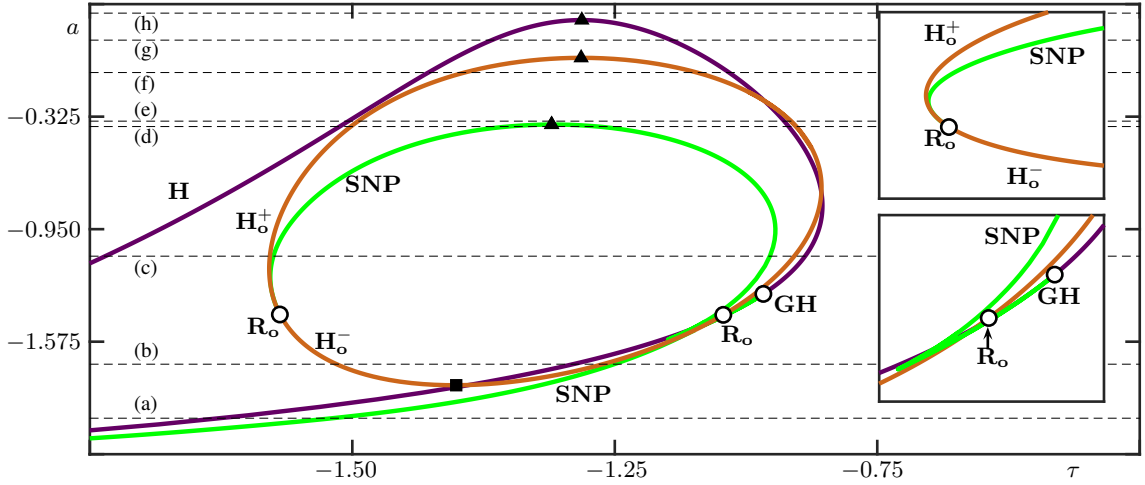}
\caption{Bifurcation diagrams in the $(\tau, a)$-parameter plane of system~\eref{eq:sanDelay} near orientable resonant bifurcation points $\mathbf{R_o}$.  Shown are the curve of homoclinic bifurcations~$\mathbf{H_o}$(brown), saddle-node bifurcations of periodic orbits $\mathbf{SNP}$  (green), and Hopf bifurcations $\mathbf{H}$ (purple).  The right superscript on $\mathbf{H_o}$ indicates the sign of the saddle quantity. The insets show enlargements close to the resonant points $\mathbf{R_o}$ (white dots) and the generalized Hopf bifurcation point $\mathbf{GH}$ (also shown as a white dot).  The black square indicates the minimum point $\mathbf{m[H_o^-]}$ in $a$ of $\mathbf{H_o^-}$, while black triangles indicate the maxima $\mathbf{M[H]}$, $\mathbf{M[SNP]}$ and $\mathbf{M[H^+_o]}$ in $a$ of $\mathbf{H}$, $\mathbf{SNP}$ and $\mathbf{H^+_o}$,  respectively.   The dashed curves indicate the slices taken for \fref{fig:bifDelayOriResSlice1} and \fref{fig:bifDelayOriResSlice2}. Here, $(\mu,\tilde{\mu},\kappa) \approx (-0.2,-0.550,-1.0)$ and other parameters are given in \Sref{sec:NoDelaySan}.} \label{fig:bifDelayOriRes} 
\end{figure} 
%%%%%%%%%%%%%%%%%%%%%%%%%%%%%%%%%%%%%%%%%%%%%%%%%%%%%%%%

\Fref{fig:bifDelayOriRes} shows the bifurcation diagram of the homoclinic bifurcation curve (brown), saddle-node bifurcation curves of periodic orbits $\mathbf{SNP}$ (green) and the Hopf bifurcation curve $\mathbf{H}$ (purple) in the $(\tau, a)$-parameter plane. Most prominently, the homoclinic orbit curve closes up on itself. As per construction, the saddle quantity $\sigma$ is independent of $\tau$ such that we encounter resonant points $\mathbf{R_o}$ at $a=a_R \approx -1.5006$, which separate the homoclinic bifurcation curve into segments $\mathbf{H_o^+}$ and $\mathbf{H_o^-}$ reflecting the sign of $\sigma$. Similar to the ODE case, a curve $\mathbf{SNP}$ emerges from each point $\mathbf{R_o}$, thus recovering the generic unfolding of an orientable resonant bifurcation locally near the points $\mathbf{R_o}$ in system~\eref{eq:sanDelay}; see the insets in \fref{fig:bifDelayOriRes}. Notice the existence of a Hopf bifurcation curve $\mathbf{H}$ of an equilibrium co-existing with $\widehat{\mathbf{0}}$. This Hopf bifurcation curve exhibits a codimension-two generalized Hopf bifurcation $\mathbf{GH}$ that changes the criticality of $\mathbf{H}$ and causes one of the $\mathbf{SNP}$ curves to coalesce (the one emanating from the left-most point $\mathbf{R_o}$). 

\subsubsection*{One-parameter bifurcation diagrams in $\tau$ for fixed $a$}

In what follows, we investigate how the orientable resonant bifurcation points organize the existence of infinitely many periodic solutions in $\tau$. As such, we choose representative values of $a$, highlighted as dashed curves in \fref{fig:bifDelayOriRes}, and present one-parameter bifurcation diagrams in $\tau$ that give qualitatively different bifurcation scenarios.

%%%%%%%%%%%%%%%%%%%%%%%%%%%%%%%%%%%%%%%%%%%%%%%%%%%%%%%%
\begin{figure}
\centering
\includegraphics{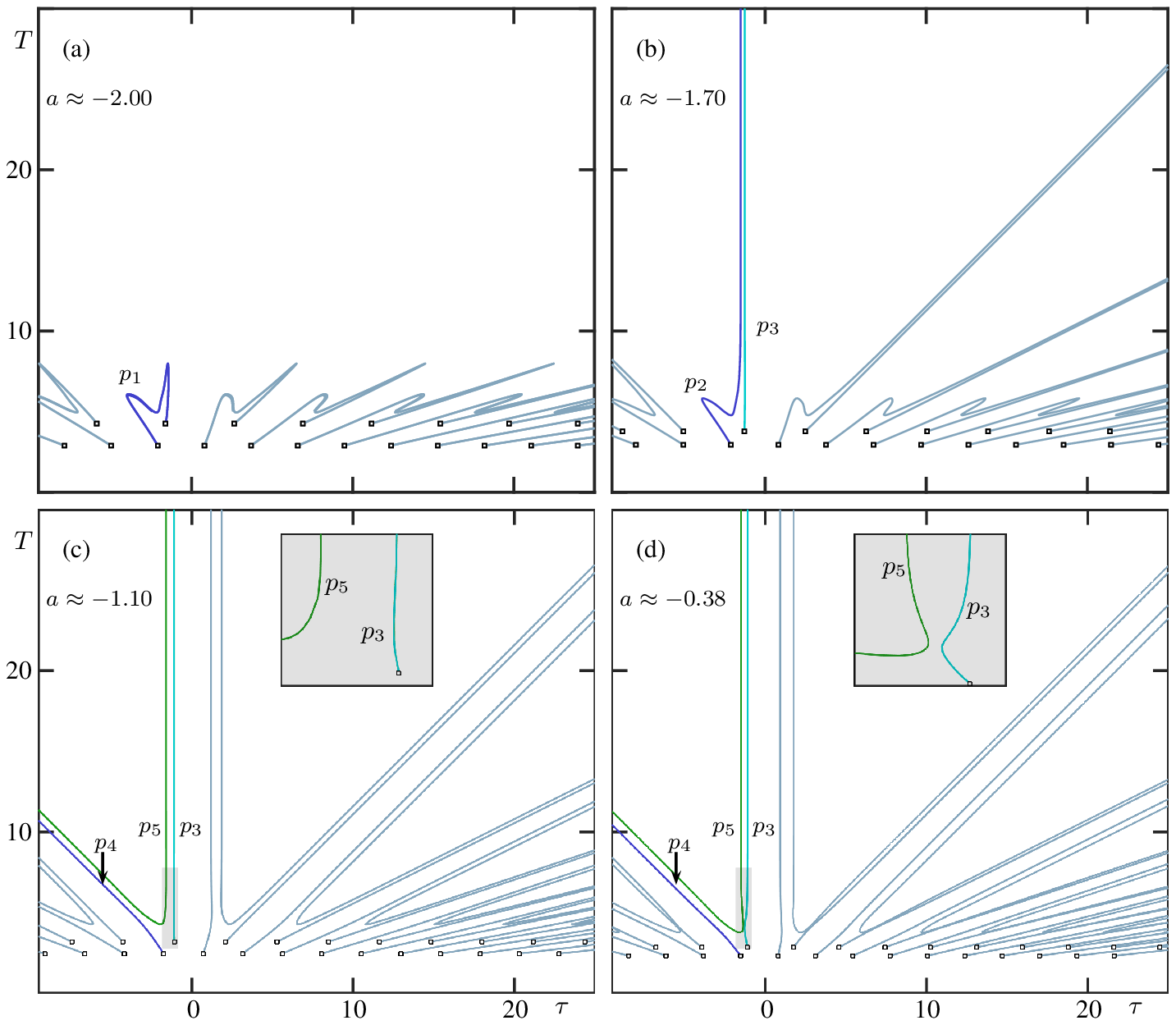}
\caption{One-parameter bifurcation diagrams in $\tau$ of system~\eref{eq:sanDelay} for different values of $a$ as indicated by the dashed lines in \fref{fig:bifDelayOriRes}. The value of $a$ for each slice is indicated on the top left side of each panel. Shown is the period $T$ of the periodic orbit with respect to $\tau$.  The branch of periodic solutions emanating from the Hopf bifurcation $\mathbf{H}$ and the homoclinic bifurcation~$\mathbf{H_o}$ shown in  \fref{fig:bifDelayOriRes} are colored blue, green and cyan, while the other periodic branches are colored gray.  The Hopf bifurcation points are indicated by small squares. } \label{fig:bifDelayOriResSlice1} 
\end{figure} 
%%%%%%%%%%%%%%%%%%%%%%%%%%%%%%%%%%%%%%%%%%%%%%%%%%%%%%%%

%%%%%%%%%%%%%%%%%%%%%%%%%%%%%%%%%%%%%%%%%%%%%%%%%%%%%%%%
\begin{figure}
\centering
\includegraphics{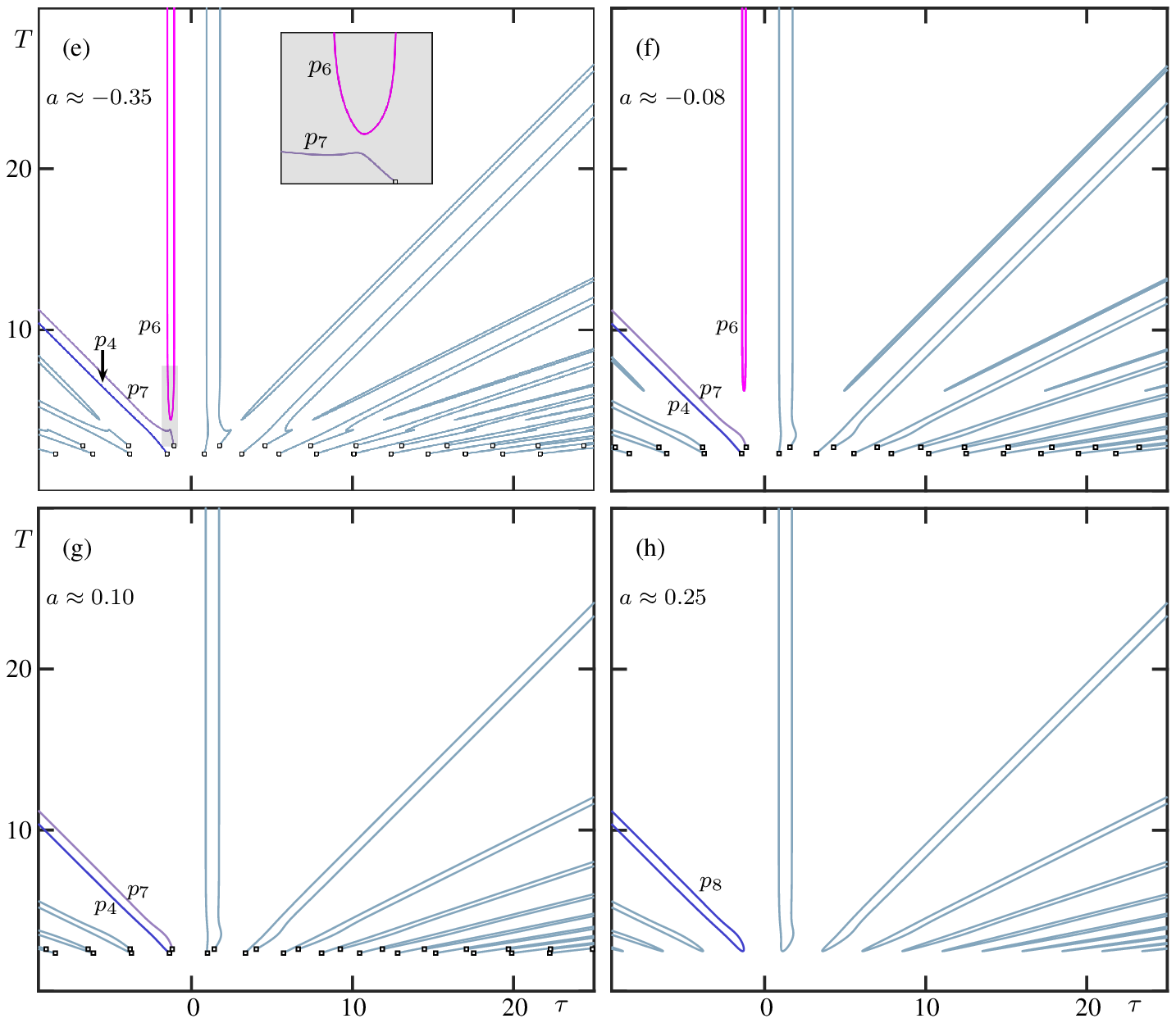}
\caption{Continued from \fref{fig:bifDelayOriResSlice1}.   The branch of periodic solutions emanating from the Hopf bifurcation $\mathbf{H}$ and the homoclinic bifurcation~$\mathbf{H_o}$ shown in  \fref{fig:bifDelayOriRes} are colored blue, magenta and light-purple,  while the other periodic branches are colored gray. } \label{fig:bifDelayOriResSlice2} 
\end{figure} 
%%%%%%%%%%%%%%%%%%%%%%%%%%%%%%%%%%%%%%%%%%%%%%%%%%%%%%%%

\Fref{fig:bifDelayOriResSlice1}(a) shows the one-parameter bifurcation diagram in $\tau$ for $a=-2$, which corresponds to the dashed line labelled (a) in \fref{fig:bifDelayOriRes}.  In \Fref{fig:bifDelayOriResSlice1}(a), system~\eref{eq:sanDelay} exhibits two Hopf bifurcations $\mathbf{H}$ at $\tau = \tau_1 \approx -2.1027$ and $\tau = \tau_2 \approx -1.6226$, which are connected through a family of periodic solutions $p_1$ (blue). The Hopf bifurcations reappear at delay values $\tau = \tau_1+2\pi k/\omega_1$ and $\tau = \tau_2+2\pi k/\omega_2$, where $k \in \Z$, and $\omega_1$ and $\omega_2$ are the imaginary parts of the respective eigenvalues at the Hopf bifurcation.  Similarly, the family of periodic orbits $p_1$ reappears at delay values $\tau =\tau_0+ k T$ for $k \in \Z$ (grey curves), where  $T$ is the period of the respective periodic orbit at $\tau=\tau_0$. However, we remark that $p_1$ can only be locally parametrized by $\tau$, since there exist values of $\tau$ where multiple periodic orbits coexist. Notice here that along the family $p_1$ and all its reappearing branches, we encounter a periodic orbit with maximum period $T_{\rm max}$. Furthermore, for sufficiently large $k$, we find saddle-node bifurcations for the reappearing families close nearby the period maxima $T_{\rm max}$ as expected from theory \cite{Yanchuk2009}.

Now we focus our attention on the qualitative changes of the one-parameter bifurcation diagram as $a$ changes. As $a$ increases,  system~\eref{eq:sanDelay} exhibits a first homoclinic bifurcation $\mathbf{m[H^-_{o}]}$ at the minimum $a$-value along the curve $\mathbf{H^-_o}$ in \fref{fig:bifDelayOriRes}. \Fref{fig:bifDelayOriResSlice1}(b) shows the situation for $a\approx -1.70$ just after crossing $\mathbf{m[H^-_{o}]}$. Here, the family of periodic solution $p_1$ split into two distinct families $p_2$ and $p_3$, each arising from a Hopf bifurcation. Notice that for both families, there does not exist a maximum period, as each family approaches a homoclinic bifurcation $\mathbf{H^-_o}$ as indicated by the dashed line~(b) in \fref{fig:bifDelayOriRes}.  The existence of $p_2$ and $p_3$ can be intuitively understood by following $T_{\rm max}$ along $p_1$. Indeed,  as $a$ increases from $-2$, the maximum period $T_{\rm max}$ growths without bounds as  system~\eref{eq:sanDelay} exhibits $\mathbf{m[H^-_{o}]}$ at $a\approx-1.8168$. Effectively, at this moment, the periodic solution associated with $T_{\rm max}$ becomes a homoclinic orbit. By increasing $a$ past $\mathbf{m[H^-_{o}]}$, the family $p_1$ splits into the families $p_2$ and $p_3$, and $\mathbf{m[H^-_{o}]}$ creates two homoclinic bifurcations $\mathbf{H^-_{o}}$. Indeed, $\mathbf{m[H^-_{o}]}$ corresponds to a fold of homoclinic bifurcations in the extended system~\eref{eq:sanDelayExt}. This transition has several important implications; most prominently, the families of periodic orbits that reappear from $p_2$ and $p_3$ in \fref{fig:bifDelayOriResSlice1}(b) give rise to temporal dissipative solitons (TDSs). For example, by applying the reappearance rule once to $p_2$,  we obtain a family of periodic orbits at $\tau=\tau_0+T$ along which the period scales linearly with $\tau$ and, therefore, the period grows beyond bound as $\tau\to\infty$. As such, the fold of homoclinic bifurcations $\mathbf{m[H^-_{o}]}$ of the extended system is a \textbf{bifurcation that gives rise to TDS}, and coincides with the bifurcation of countably-many saddle-node bifurcations with infinite period. Effectively, as we near the point $\mathbf{m[H^-_{o}]}$, the saddle-node bifurcations shown in panel~(a) for each reappearing family occur at larger and larger values of the delay until they disappear at infinity thus producing TDSs.

Notice also that the first reappearance of family $p_2$ attains a local period-maximum for small positive $\tau$ as observed in \fref{fig:bifDelayOriResSlice1}(b). Analogous to the transition between panels~(a) and (b), this local maximum increases without bounds as $a$ increases; thus, implying the existence of another fold of homoclinic bifurcations in the extended system~\eref{eq:sanDelayExt} for small positive $\tau$ (not shown in \fref{fig:bifDelayOriRes}) that must have occurred as we transition between panels~(b) and (c). Such fold splits the family $p_2$ into families $p_4$ and $p_5$. As before, we see that this transition is accompanied by countably-many saddle-node bifurcations with infinite period; particularly, the saddle-node bifurcation of $p_2$. At this secondary fold of homoclinic bifurcations in the extended system~\eref{eq:sanDelayExt}, $p_2$ extends towards minus infinity and splits into $p_4$ and $p_5$. Panel~(c) corresponds to the situation past the resonant point $\mathbf{R_o}$ at $a \approx -1.1$. This implies the existence of saddle-node periodic orbits that exist in families $p_3$ and also $p_5$. However, as these saddle-node periodic orbits emanate from the resonant points $\mathbf{R_o}$ with infinite period, we cannot observed them in panel~(c). To visualise the existence of these saddle-node periodic orbits, we increase $a$ to $a \approx -0.38$ and show in panel~(d) the situation farther away from point $\mathbf{R_o}$. Here, these saddle-node periodic orbits along the families $p_3$ and $p_5$ are clearly visible. Notice also that these folds occur very close to each other in the $\tau$ parameter.

Recall from \fref{fig:bifDelayOriRes} that the saddle-node bifurcation curve $\mathbf{SNP}$ reaches a maximum $\mathbf{M[SNP]}$ in $a$, where $a\approx -0.3689$. Indeed, the situation showcased in panel~(d) occurs just before reaching this maximum. At this maximum, the saddle-node periodic orbits shown in panel~(d) coincide, and families $p_3$ and $p_5$ undergo a saddle transition. We showcase the consequence past this transition in \fref{fig:bifDelayOriResSlice2}(e) for $a \approx -0.35$; here, we clearly observe that segments of the families $p_3$ and $p_5$ have rearranged to becomes families $p_6$ and $p_7$. Of particular interest is family $p_6$, which now connects the two homoclinic bifurcations $\mathbf{H^+_o}$ for small negative values of $\tau$. Notice also that the branches reappearing from $p_6$ are not connected to any of the Hopf bifurcations. On the other hand, family $p_7$ does not exhibit a homoclinic bifurcation; however, its period scales linearly with negative values of $\tau$. Thus, its first reappearance for positive $\tau$ connects a Hopf bifurcation and a homoclinic bifurcation for small positive $\tau$. Even though the saddle transition is a local phenomenon, it effectively rearranges the families of periodic orbits in $\tau$ at a global scale.

Notice that along the family $p_6$ we encounter a minimum in $T$.  As we increase $a$, we observe that the family $p_6$ exhibits  its minimum at a larger value of $T$, see panel~(g). This is relevant because it creates a mechanism that allows for families of periodic orbits to disappear in analogy to the transition observed in panels~(a) and (b). Indeed, as we increase $a$ to the value where the curve $\mathbf{H^+_o}$ reaches its maximum, $\mathbf{M[H^+_o]}$, we encounter another fold of homoclinic bifurcations in the extended system~\eref{eq:sanDelayExt}. As a result, the family $p_6$ vanishes with infinite period, and it is absent in panel~(g) where we have passed this maximum. For completeness, we showcase the situation past the maximum of the Hopf bifurcation curve~$\mathbf{M[H]}$, see panel~(h). Past this maximum, the Hopf bifurcation points have disappeared  from the bifurcation diagram pairwise, and the family $p_4$ and $p_7$ have merged into a new family $p_8$.  Notice that family $p_8$ now has a saddle-node bifurcation (not indicated) which is mapped close to the minimum of its first reappearing branch, which reappears infinitely-many times throughout the bifurcation diagram.

To summarize, although $\mathbf{m[H^-_o]}$  and $\mathbf{M[H^+_o]}$ both correspond to folds of homoclinic bifurcations in the extended system~\eref{eq:sanDelayExt}, they mediate different transitions that organize families of periodic solutions. Indeed, the transition observed through the minimum $\mathbf{m[H^-_o]}$ takes a family of periodic orbits with finite period and splits it into two families, which are disconnected by two homoclinic bifurcations. In this sense, the minimum $\mathbf{m[H^-_o]}$ creates/destroys (depending on the direction of $a$) a $\tau$-gap in a family of periodic orbits for a fixed $a$. On the other hand, the transition observed through the maximum $\mathbf{M[H^+_o]}$ creates/destroys a family of periodic orbits that connects two homoclinic bifurcations. Loosely speaking, one transition creates/destroys a gap in a family of periodic orbits, and the other creates/destroys a family of periodic orbits.

%%%%%%%%%%%%%%%%%%%%%%%%%%%%%%%%%%%%%%%%%%%%%%%%%%%%%%%%
\begin{table}[]
\centering
\begin{tabular}{|ccccccccl|}
\hline
\multicolumn{5}{|c|}{\textbf{Parameters}} & \multicolumn{1}{c|}{\multirow{2}{*}{\textbf{\begin{tabular}[c]{@{}c@{}}Continuation\\ parameters\end{tabular}}}} & \multicolumn{3}{c|}{\multirow{2}{*}{\textbf{\begin{tabular}[c]{@{}c@{}}Stopping \\ condition\end{tabular}}}} \\ \cline{1-5}
\multicolumn{1}{|c|}{$a$} & \multicolumn{1}{c|}{$\mu$} & \multicolumn{1}{c|}{$\tilde{\mu}$} & \multicolumn{1}{c|}{$\kappa$} & \multicolumn{1}{c|}{$\tau$} & \multicolumn{1}{c|}{} & \multicolumn{3}{c|}{} \\ \hline
\multicolumn{9}{|c|}{\textbf{Homoclinic solution}} \\ \hline
\multicolumn{1}{|c|}{$-2.0620$} & \multicolumn{1}{c|}{$0.0550$} & \multicolumn{1}{c|}{$0.03$} & \multicolumn{1}{c|}{$0$} & \multicolumn{1}{c|}{$0$} & \multicolumn{1}{c|}{$\kappa,\mu$} & \multicolumn{1}{r}{$\mu$} & $=$ & $0.2$ \\ \hline
\multicolumn{1}{|c|}{$-2.0620$} & \multicolumn{1}{c|}{$0.2$} & \multicolumn{1}{c|}{$0.03$} & \multicolumn{1}{c|}{$-0.3805$} & \multicolumn{1}{c|}{$0$} & \multicolumn{1}{c|}{$\kappa,\tau$} & \multicolumn{1}{r}{$\kappa$} & $=$ & $-1$ \\ \hline
\multicolumn{1}{|c|}{\multirow{2}{*}{$-2.0620$}} & \multicolumn{1}{c|}{\multirow{2}{*}{$0.2$}} & \multicolumn{1}{c|}{\multirow{2}{*}{$0.03$}} & \multicolumn{1}{c|}{\multirow{2}{*}{$-1.0$}} & \multicolumn{1}{c|}{$-0.9114$} & \multicolumn{1}{c|}{\textemdash} & \multicolumn{3}{c|}{\textemdash} \\ \cline{5-9} 
\multicolumn{1}{|c|}{} & \multicolumn{1}{c|}{} & \multicolumn{1}{c|}{} & \multicolumn{1}{c|}{} & \multicolumn{1}{c|}{$0.9031$} & \multicolumn{1}{c|}{\textemdash} & \multicolumn{3}{c|}{\textemdash} \\ \hline
\multicolumn{9}{|c|}{\textbf{Two-homoclinic solution}} \\ \hline
\multicolumn{1}{|c|}{$-0.7490$} & \multicolumn{1}{c|}{$0.0550$} & \multicolumn{1}{c|}{$0.03$} & \multicolumn{1}{c|}{$0$} & \multicolumn{1}{c|}{$0$} & \multicolumn{1}{c|}{$\kappa,\mu$} & \multicolumn{1}{r}{$\mu$} & $=$ & $0.2$ \\ \hline
\multicolumn{1}{|c|}{$-0.7490$} & \multicolumn{1}{c|}{$0.2$} & \multicolumn{1}{c|}{$0.03$} & \multicolumn{1}{c|}{$-0.5613$} & \multicolumn{1}{c|}{} & \multicolumn{1}{c|}{$\kappa,\tau$} & \multicolumn{1}{r}{$\kappa$} & $=$ & $-1$ \\ \hline
\multicolumn{1}{|c|}{\multirow{2}{*}{$-0.7490$}} & \multicolumn{1}{c|}{\multirow{2}{*}{$0.2$}} & \multicolumn{1}{c|}{\multirow{2}{*}{$0.03$}} & \multicolumn{1}{c|}{\multirow{2}{*}{$-1.0$}} & \multicolumn{1}{c|}{$-0.4384$} & \multicolumn{1}{c|}{\textemdash} & \multicolumn{3}{c|}{\textemdash} \\ \cline{5-9} 
\multicolumn{1}{|c|}{} & \multicolumn{1}{c|}{} & \multicolumn{1}{c|}{} & \multicolumn{1}{c|}{} & \multicolumn{1}{c|}{$0.4309$} & \multicolumn{1}{c|}{\textemdash} & \multicolumn{3}{c|}{\textemdash} \\ \hline
\multicolumn{9}{|c|}{\textbf{Period-doubling bifurcation}} \\ \hline
\multicolumn{1}{|c|}{$-0.6414$} & \multicolumn{1}{c|}{$0.0550$} & \multicolumn{1}{c|}{$0.03$} & \multicolumn{1}{c|}{$0$} & \multicolumn{1}{c|}{$0$} & \multicolumn{1}{c|}{$\kappa,\mu$} & \multicolumn{1}{r}{$\mu$} & $=$ & $0.2$ \\ \hline
\multicolumn{1}{|c|}{$-0.6414$} & \multicolumn{1}{c|}{$0.2$} & \multicolumn{1}{c|}{$0.03$} & \multicolumn{1}{c|}{$-0.5679$} & \multicolumn{1}{c|}{$0$} & \multicolumn{1}{c|}{$\kappa,\tau$} & \multicolumn{1}{r}{$\kappa$} & $=$ & $-1$ \\ \hline
\multicolumn{1}{|c|}{\multirow{2}{*}{$-0.6414$}} & \multicolumn{1}{c|}{\multirow{2}{*}{$0.2$}} & \multicolumn{1}{c|}{\multirow{2}{*}{$0.03$}} & \multicolumn{1}{c|}{\multirow{2}{*}{$-1.0$}} & \multicolumn{1}{c|}{$-0.4250$} & \multicolumn{1}{c|}{\textemdash} & \multicolumn{3}{c|}{\textemdash} \\ \cline{5-9} 
\multicolumn{1}{|c|}{} & \multicolumn{1}{c|}{} & \multicolumn{1}{c|}{} & \multicolumn{1}{c|}{} & \multicolumn{1}{c|}{$0.4223$} & \multicolumn{1}{c|}{\textemdash} & \multicolumn{3}{c|}{\textemdash} \\ \hline
\end{tabular}\vspace{6pt}
\caption{Additional information for the computation of the bifurcation diagram for the non-orientable resonant case in system~\eref{eq:sanDelay}.  Shown are the sequences of continuation runs to lift the homoclinic bifurcation and bifurcations of periodic solutions to non-zero $\tau$. Other parameters are given in \Sref{sec:NoDelaySan}.}\label{tab:nonOriSeqCont}
\end{table}

%%%%%%%%%%%%%%%%%%%%%%%%%%%%%%%%%%%%%%%%%%%%%%%%%%%%%%%%

\subsection{Codimension-two non-orientable resonant case}\label{sec:NonOrientableResonantCase}

\begin{figure}
\centering
\includegraphics{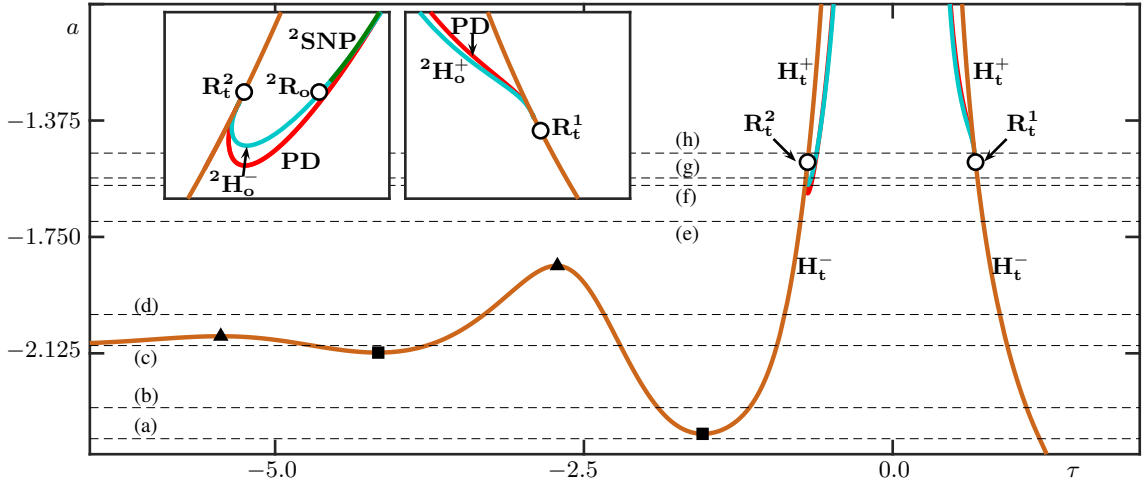}
\caption{Bifurcation diagrams in the $(\tau, a)$-parameter plane of system~\eref{eq:sanDelay} near two non-orientable resonant homoclinic bifurcations $\mathbf{R^1_t}$ and $\mathbf{R^2_t}$.  Shown are the curve of homoclinic bifurcations~$\mathbf{H_t}$(brown), period-doubling bifurcations $\mathbf{PD}$ (red), and two-homoclinic bifurcations $\mathbf{H^2_o}$ (cyan).  The inset show magnifications close to the resonant points $\mathbf{R^1_t}$ and $\mathbf{R^2_t}$ which are represented by white dots.  The black square indicates the local minima $\mathbf{m[H_o^-]}$ in $a$ of $\mathbf{H_t^-}$, while the black triangles indicate maxima $\mathbf{M[H^-_t]}$.  The dashed curves indicate the slices considered in \fref{fig:bifDelayNonOriResSlice1} and \fref{fig:bifDelayNonOriResSlice2}. Here, $(\mu,\tilde{\mu},\kappa) \approx (0.2,0.03,-1.0)$ and other parameters are given in \Sref{sec:NoDelaySan}.} \label{fig:bifDelayNonOriRes} 
\end{figure}

Here, we present the numerical bifurcation diagram of the non-orientable homoclinic bifurcation in the $(\tau,a)$-parameter plane for fixed parameters as in \cref{sec:OrientableResonantCase}. Our procedure is in part analogous to the orientable resonant case; that is, we first choose a homoclinic orbit of system \cref{eq:san} at parameter values 
$$(a,\mu,\tilde \mu)=(-2.0620,0.0550,0.03)$$ computed using \textsc{Auto7p} in system \cref{eq:san} and then lift the homoclinic orbit to system \cref{eq:sanDelay} in a series of two-parameter continuation runs in DDE-BIFTOOL as a periodic orbit with fixed, large period. In particular, we first `switch on' $\kappa$ by continuing in $\kappa$ and $\mu$, and then $\tau$ by continuing in $\kappa$ and $\tau$; consider \cref{tab:nonOriSeqCont} for details. We observe that the curve of homoclinic orbits in $\kappa$ and $\tau$ folds back on itself, so that we encounter two different homoclinic bifurcation points at the chosen stopping condition $\kappa=-1$. Namely, at parameter values $(a,\mu,\tilde\mu,\kappa,\tau)$ approximately equal to 
$$(-2.0620,0.3,0.03,-1.0,-0.9114)\quad\mbox{and}\quad 
(-2.0620,0.3,0.03,-1.0,0.9031).$$
Unlike the previous section, we obtain two different starting homoclinic orbits that give rise to two different curves of homoclinic bifurcations in the $(\tau,a)$-parameter plane. \Fref{fig:bifDelayNonOriRes} shows the resulting two curves (both labeled) $\mathbf{H_t}$ (brown) of homoclinic bifurcations, which are not connected in $(\tau,a)$-parameter plane, and they exist either for negative or positive delays. Similar to before, the superscript of $\mathbf{H_t}$ indicates the sign of the saddle quantity of the origin. It follows from Eq.~\cref{eq:EigCond} that the saddle quantity becomes zero at $a = a_R \approx -1.5007$; thus a codimension-two non-orientable resonant bifurcation point $\mathbf{R_t}$ occurs along each of the curves $\mathbf{H_t}$. We use a superscript to distinguish these two points. From each of the points $\mathbf{R_t},$ we find a curve of period-doubling bifurcations $\mathbf{PD}$ and a curve of orientable two-homoclinic bifurcation $\mathbf{^2H_o}$; see the inset of \fref{fig:bifDelayNonOriRes}.  Again, we use a superscript on the homoclinic bifurcation curves to indicate the sign of the saddle quantity.

Notice that the situation at the two points $\mathbf{R_t}$ corresponds to different unfoldings of the non-orientable resonant homoclinic bifurcation. While the curves $\mathbf{PD}$ and $\mathbf{^2H_o^+}$ emerge with positive saddle quantity from the point $\mathbf{R_t^1}$, the bifurcation curves emerging from $\mathbf{R_t^2}$ emanate with negative saddle quantity. This local difference can be understood as a consequence of higher order terms in the Poincar\'{e} map as in the ODE case \cite{Chow1990}. Away from $\mathbf{R_t^2}$, we observe that the curves $\mathbf{PD}$ and $\mathbf{^2H_o^+}$ fold back on themselves in the $(\tau,a)$-parameter plane. This implies the existence of a secondary resonance point $\mathbf{^2R_o}$ when the curve $\mathbf{^2H_o^+}$ passes through the resonance condition $a=a_R$. Since points along $\mathbf{^2H_o^+}$ correspond to orientable homoclinic bifurcations, a saddle-node bifurcation curve $\mathbf{^2SNP}$ (dark-green) of period-doubled periodic orbits emanates from $\mathbf{^2R_o}$, see the left inset in \fref{fig:bifDelayNonOriRes}.

%%%%%%%%%%%%%%%%%%%%%%%%%%%%%%%%%%%%%%%%%%%%%%%%%%%%%%%% 
%%%%%%%%%%%%%%%%%%%%%%%%%%%%%%%%%%%%%%%%%%%%%%%%%%%%%%%%

For the parameter range chosen in \fref{fig:bifDelayNonOriRes}, notice that the left curve $\mathbf{H_t^-}$ extends to large negative values of $\tau$ in an oscillatory fashion and asymptotes at $a\approx  -2.09$.
This structure has several implications for the existence of periodic orbits with large period near the curve $\mathbf{H_t^-}$. Like it was shown in \cref{fig:bifDelayOriResSlice1,fig:bifDelayOriResSlice2} for the resonant case, each minimum $\mathbf{m}[\mathbf{H_t^-}]$ and maximum $\mathbf{M}[\mathbf{H_t^-}]$ in $a$ along the curve $\mathbf{H_t^-}$ corresponds to a fold of homoclinic bifurcations in the extended system \cref{eq:sanDelayExt}.

%%%%%%%%%%%%%%%%%%%%%%%%%%%%%%%%%%%%%%%%%%%%%%%%%%%%%%%%
\begin{figure}
\centering
\includegraphics{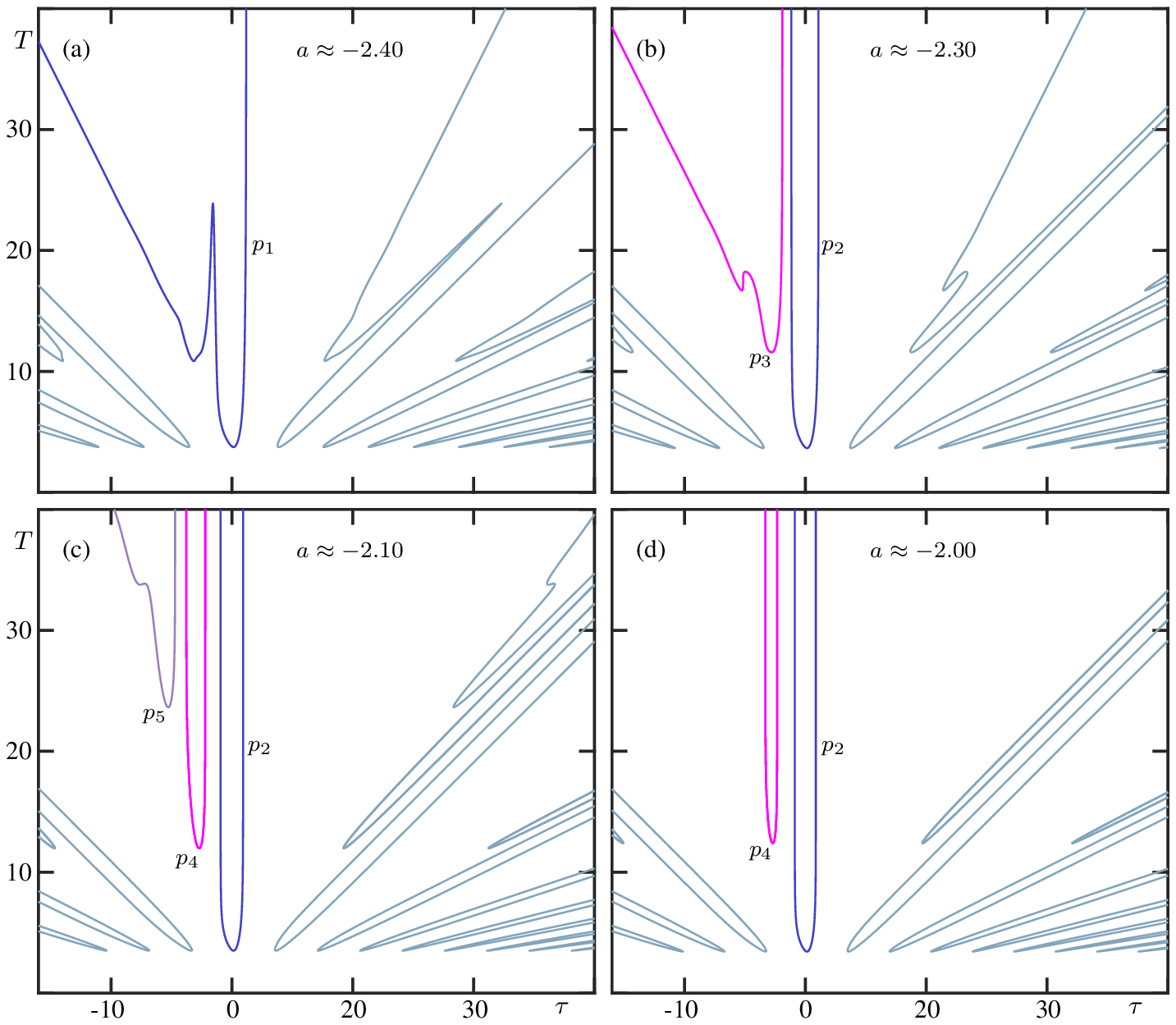}
\caption{One-parameter bifurcation diagrams in $\tau$ of system~\eref{eq:sanDelay} for different values of $a$ as indicated in \fref{fig:bifDelayNonOriRes}. The value of $a$ for each slice is indicated on the top left side of each panel. Shown is the period $T$ of the periodic orbit with respect to $\tau$.   The branch of periodic solutions emanating from the homoclinic bifurcations~$\mathbf{H_t}$ are colored blue, magenta and light-purple, while the other periodic branches are colored gray.}\label{fig:bifDelayNonOriResSlice1} 
\end{figure} 
%%%%%%%%%%%%%%%%%%%%%%%%%%%%%%%%%%%%%%%%%%%%%%%%%%%%%%%%

%%%%%%%%%%%%%%%%%%%%%%%%%%%%%%%%%%%%%%%%%%%%%%%%%%%%%%%%
\begin{figure}
\centering
\includegraphics{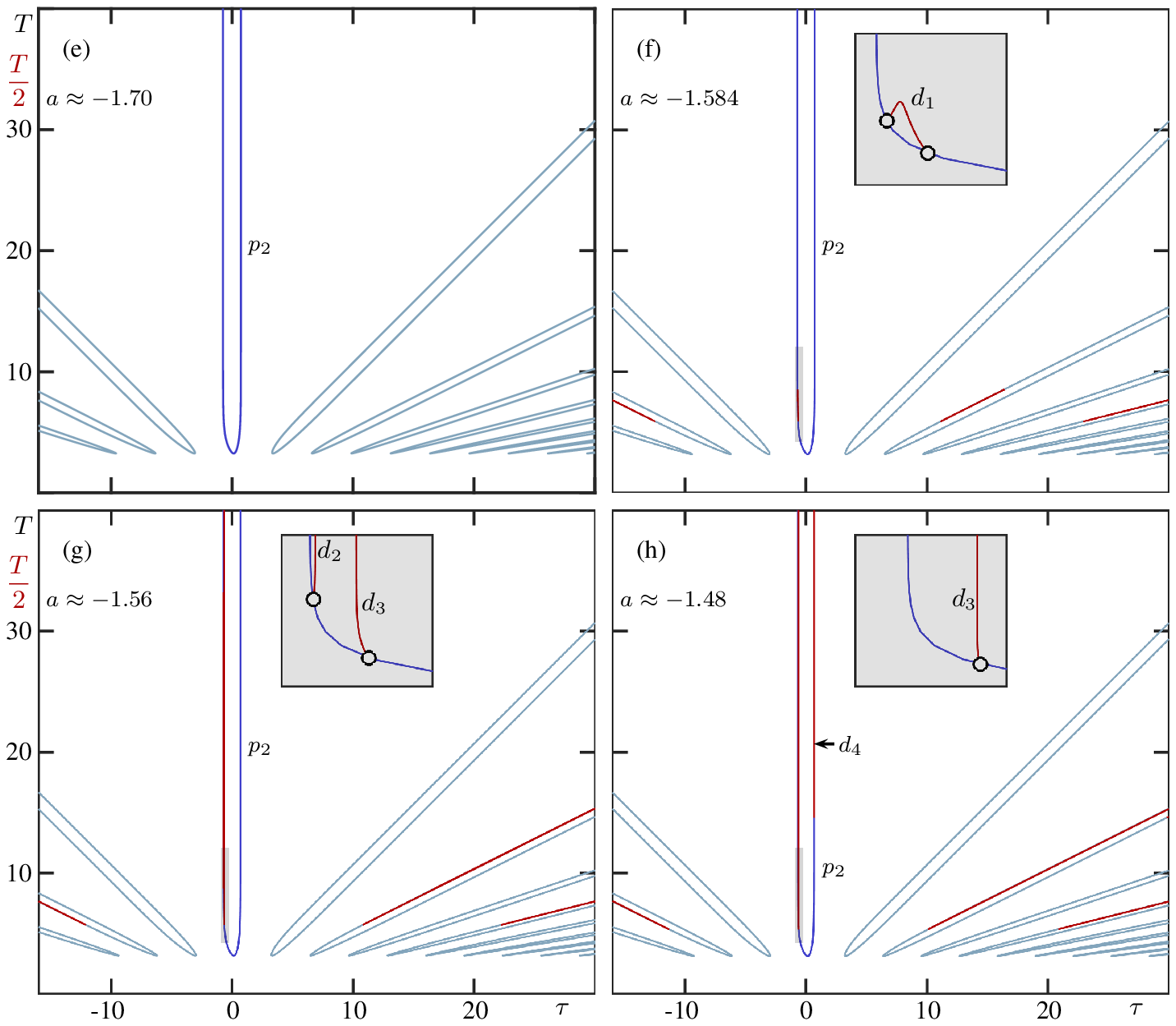}
\caption{Continued from \fref{fig:bifDelayNonOriResSlice1}. The branch of periodic solutions emanating from the homoclinic bifurcations~$\mathbf{H_t}$ shown in  \fref{fig:bifDelayOriRes} is colored blue, and the ones emanating from the two-loop homoclinic bifurcation $\mathbf{^2H^+_o}$ and period-doubling bifurcation $\mathbf{PD}$ by red, while the other periodic branches are colored gray. }\label{fig:bifDelayNonOriResSlice2} 
\end{figure} 
%%%%%%%%%%%%%%%%%%%%%%%%%%%%%%%%%%%%%%%%%%%%%%%%%%%%%%%%

We now study the existence of periodic orbits near and far away from the non-orientable resonance points. We take representative slices for fixed values of $a$ as indicated by the dashed lines in \fref{fig:bifDelayNonOriRes} as in \cref{sec:OrientableResonantCase} for the orientable resonant case. \Fref{fig:bifDelayNonOriResSlice1}(a) shows the $\tau$-bifurcation diagram for $a \approx -2.40$. Notice that this slice intersects $\mathbf{H^-_t}$ only once and, thus, a family of periodic orbits denoted by $p_1$ (blue curve) emerges from the point of bifurcation. Notice also in panel~(a) that we observe two local period minima and one local period maximum as $\tau$ decreases along $p_1$. Interestingly, the period of the family $p_1$ scales as $T/|\tau| \approx 2$ as $\tau$ decreases towards minus infinity, see \fref{fig:bifInfCycle}(a) for the profile of the periodic orbit along $p_1$ at $\tau \approx -100$ and \cref{sec:twoInfty} for discussion of this limit. Panel~(a) also shows the branches reappearing from $p_1$ at other values of the delay (gray). It is worth noting here that the first appearance for positive $\tau$ of $p_1$ satisfies two different scalings as $\tau$ goes to infinity: $T/\tau \approx 1$ which corresponds to the reappearance of the periodic orbits close to the homoclinic in $p_1$, and $T/ \tau \approx 2$. Also, similar to our discussion of the orientable case, the local period extrema along $p_1$ reappear close to points near saddle-node bifurcations \cite{Yanchuk2009}.

As we increase $a$, the family $p_1$ exhibits a fold of homoclinic bifurcations $\mathbf{m[H^-_t]}$ in the extended system~\eref{eq:sanDelayExt}.  \Fref{fig:bifDelayNonOriResSlice1}(b) shows the $\tau$ bifurcation diagram at $a\approx -2.30$, that is, past this bifurcation point. We observe in this panel that the local maximum in panel~(a) has increased without bounds, and the family $p_1$ has split in two families $p_2$ (blue) and $p_3$ (magenta) which are separated by means of two homoclinic bifurcations created at $\mathbf{m[H^-_t]}$. As a result, the family $p_2$ is bounded in $\tau$ by two homoclinic bifurcations, while the family $p_3$ is of the same form as family $p_1$ in panel~(a). This new configuration carries over to the reappearing branches of $p_2$ and $p_3$; indeed, the reappearing families are now disconnected and their first reappearances scale as either $T/\tau \approx 1$ or $T/\tau \approx 2$ when $\tau$ goes to infinity. Remarkably, this process of branch splitting continues as we further increase $a$. Indeed, the curve $\mathbf{H^-_t}$ exhibits a second local minimum $\mathbf{m[H^-_t]}$ in $a$, see \fref{fig:bifDelayNonOriRes}, which corresponds to a second fold of homoclinic bifurcations $\mathbf{m[H^-_t]}$ in the extended system~\eref{eq:sanDelayExt}. Past this point, at $a\approx -2.10$, the family $p_3$ has split into families $p_4$ and $p_5$ as can be observed in \Fref{fig:bifDelayNonOriResSlice1}(c).

Recall that the curve $\mathbf{H^-_t}$ oscillates around the asymptotic value $a\approx -2.09$ in the $(\tau,a)$-parameter plane as $\tau$ goes to minus infinity. This suggests that the branch splitting process occurs countably many times whenever we cross a minimum of the curve $\mathbf{H^-_t}$. Thus, whenever the extended system~\eref{eq:sanDelayExt} exhibits a fold of homoclinic bifurcations, a new family of periodic orbits emerges in-between the $\tau$ values where the newly created homoclinic orbits occur for fixed $a$. Our numerics suggest that at the asymptotic value $a \approx -2.09$, there does not exists a family scaling as $T / | \tau | \approx 2$ as $\tau$ goes to minus infinity, see \cref{sec:twoInfty} for more details.  Past the asymptote, we have numerical evidence of countably many maxima in $a$ occurring along curve $\mathbf{H^-_t}$. Each of these maxima causes one of the families created by crossing a minimum to vanish as $a$ increases past the maximum value, analogous to what was observed for the orientable resonant case in \cref{sec:OrientableResonantCase}. To showcase this mechanism, we present in \Fref{fig:bifDelayNonOriResSlice1}(d) at $a\approx -2.00$ the $\tau$-bifurcation diagram past the left-most maximum $\mathbf{M[H^-_t]}$ in \fref{fig:bifDelayNonOriRes}. As a result of this process, the family $p_5$ in panel~(c) has disappeared in panel~(d). Increasing $a$ past the last local maxima of $\mathbf{H^-_t}$ to $a\approx -1.70$,  we again witness the disappearance of a family of periodic orbits in \Fref{fig:bifDelayNonOriResSlice2}(e), namely, the $p_4$ family.  Between these two panels, the period minimum of the family $p_4$ grows without bounds as we approach the last local maximum $\mathbf{M[H^-_t]}$, which corresponds to a fold of homoclinic bifurcations in the extended system~\eref{eq:sanDelayExt}.

Now we continue our analysis close to the non-orientable resonants points $\mathbf{R^{1,2}_t}$ shown in \fref{fig:bifDelayNonOriRes}. Consider now \fref{fig:bifDelayNonOriResSlice2} which showcases the $\tau$-bifurcation diagrams for values of $a$ below and above the resonant points. Notice that increasing $a$ from $-1.70$ (panel~(e)) to $-1.584$ (panel~(f)), we cross a local minimum in $a$ of a curve of period-doubling bifurcations $\mathbf{PD}$ that emanates from $\mathbf{R^{2}_t}$. As a result, the $\tau$-bifurcation diagram, shown in panel~(f), exhibits two period-doubling bifurcation points which are indicated as white dots in the corresponding magnification (gray region). These bifurcation points give rise to a family $d_1$ of period-doubled orbits (red curve), where the period is represented by half of its value to emphasize the branching of $d_1$ with respect to $p_2$. Panel~(f) also shows the reappearance of family $d_1$ for different values of the delay (red). As a consequence of \cref{thm:reap-bif}, the period-doubling bifurcation points and their corresponding period-doubled families only reappear along every second branch reappearing from $p_2$. Since the family $d_1$ exhibits a local period maximum, the families reappearing from $d_1$ only exist for bounded values of $\tau$, where the maximum of $d_1$ maps near saddle-node bifurcations.

As we increase $a$, we cross  the minimum $\mathbf{m[^2H^-_o]}$  of the two-homoclinic bifurcation curve $\mathbf{^2H^-_o}$ that emanates from $\mathbf{R^{2}_t}$. At the point $\mathbf{m[^2H^-_o]}$, the extended system~\eref{eq:sanDelayExt} undergoes a fold of two-homoclinic bifurcations. \Fref{fig:bifDelayNonOriResSlice2}(g) shows the situation when $a$ is increased to $-1.56$ past $\mathbf{m[^2H^+_o]}$. As a result, family $d_1$ has split into families $d_2$ and $d_3$ separated by the two newly created two-homoclinic bifurcations in the $\tau$-bifurcation diagram, as compared to panel~(f). This transition is analogous to the already seen folds of homoclinic bifurcations $\mathbf{m[H^-_t]}$ in the extended system~\eref{eq:sanDelayExt} observed in \fref{fig:bifDelayNonOriResSlice1}. Note that the period along families $d_2$ and $d_3$ now grows beyond bound as each approaches a two-homoclinic orbit. As a direct consequence, the families reappearing from $d_2$ and $d_3$ now exist for $\tau$-values which extend towards infinity; thus, they correspond to bound two-pulse TDSs. As we transition past the $a$-value where system~\eref{eq:sanDelay} exhibits the non-orientable resonant bifurcations $\mathbf{R^{1,2}_t}$, we observe in panel~(h) that family $d_2$ has disappeared. This can be intuitively understood by the behavior of the period-doubling bifurcation curve $\mathbf{PD}$ close to $\mathbf{R^{2}_t}$, see \fref{fig:bifDelayNonOriRes}. Namely, as we approach $a_R \approx -1.5007$, the period of the orbit undergoing the  period-doubling bifurcation (left white point in the inset of panel~(g)) grows without bounds, which can be seen as the left-most period-doubling bifurcation point disappearing with infinite period in the $\tau$-bifurcation diagram at $a_R$. At the same time, a new period-doubling point emerges with infinite period along family $p_2$ for small positive values of $\tau$. Past this point, the period-doubling bifurcation creates a family $d_4$ of period-doubled orbits that connects the period-doubling bifurcation $\mathbf{PD}$ and the two-homoclinic bifurcation $\mathbf{^2H^+_o}$ as shown in panel~(h) for $a\approx -1.48$. Recall that the curves $\mathbf{PD}$ and $\mathbf{^2H^+_o}$ related to $d_4$ arise from the resonant point $\mathbf{R^1_t}$. The reappearing families of $d_4$ exist for $\tau$-values which extend towards infinity; thus, they generate another set of bound two-pulse TDSs which coexists with the ones generated by $d_3$. In this way, the non-orientable resonant bifurcation organizes the existence of bound two-pulse TDSs by creating/destroying entire families of bound two-pulse TDSs.

\subsubsection{Limiting behavior of orbits near and along $\mathbf{H^-_t}$ as $\tau \to -\infty$}\label{sec:twoInfty}
Recall that in \fref{fig:bifDelayNonOriRes} the left homoclinic curve $\mathbf{H^-_t}$ asymptotes at $a \approx -2.09$ as $\tau \rightarrow -\infty$ in an oscillatory fashion thereby creating countably many folds of homoclinic bifurcations in the extended system~\eref{eq:sanDelayExt}. To understand the underlying object that organizes this behavior, we explore in this section the homoclinic orbits along $\mathbf{H^-_t}$ for large negative values of $\tau$. We also consider the limiting behavior of periodic orbits with scaling $T / | \tau | \approx 2$ as $\tau \rightarrow -\infty$.  To achieve this, we consider the periodic family $p_1$, shown in \fref{fig:bifDelayNonOriResSlice1}(a), as a representative for such families. In this way, we can study how the limiting object of the family $p_1$, at $\tau$ minus infinity, approaches the limiting object of the homoclinic orbits along the curve $\mathbf{H^-_t}$. More precisely, we consider a representative large negative value of $\tau$, e.g., $\tau=-100$, and then study the profile of the corresponding periodic orbit of $p_1$ and the corresponding homoclinic orbit of  $\mathbf{H^-_t}$.

%Notice that if one fixes $\tau$ at a large negative value $\tau_o$ then the corresponding periodic orbit along  the family $p_1$ will converge to a homoclinic orbit  at $\tau_o$ as we increase $a$.

%%%%%%%%%%%%%%%%%%%%%%%%%%%%%%%%%%%%%%%%%%%%%%%%%%%%%%%%
\begin{figure}
\centering
\includegraphics{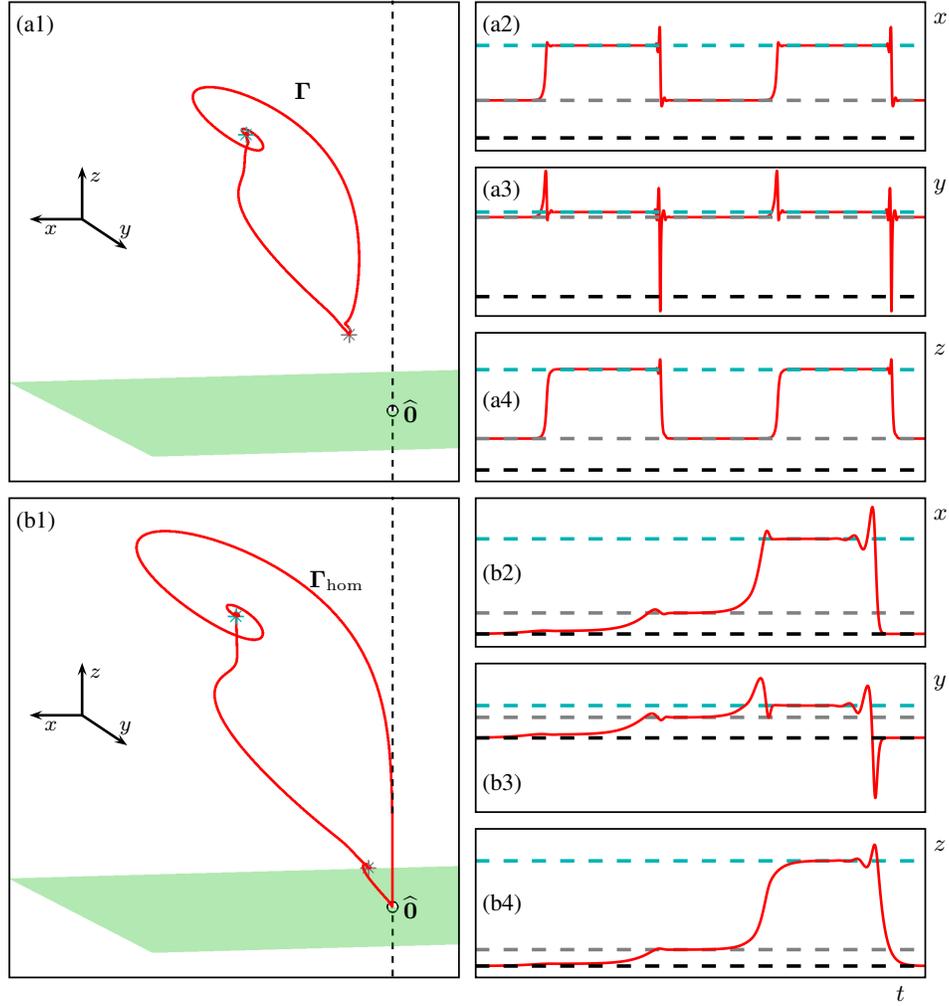}
\caption{Phase portraits  of system~\eref{eq:sanDelay} showing a periodic orbit (red curve) along the family $p_1$, panel~(a1), and a homoclinic orbit (red curve) along $\mathbf{H^+_{t}}$, panel~(b1), for a large negative value of $\tau$.  For both panels, shown are: the saddle equilibrium $\mathbf{\widehat{0}}$ (white dot), the $xy$-plane (green plane), and the computed two-period points $\mathbf{q_1}$ (cyan asterisk) and $\mathbf{q_2}$  (grey asterisk) from the multivalued map~\eref{eq:relationAtInfy}.  Panels~(a2)-(a4) and (b2)-(b4) shows the temporal traces of the periodic and homoclinic orbits, respectively. In these panels,  the origin $\mathbf{\widehat{0}}$ (black), $\mathbf{q_1}$ (cyan) and $\mathbf{q_2}$ (grey) are indicated as dashed lines.  Here,  $(a, \tau) \approx  (-2.4,  100)$ in panels~(a), and $(a, \tau) = (-2.0935,  100)$ in panels~(b).} \label{fig:bifInfCycle} 
\end{figure}
%%%%%%%%%%%%%%%%%%%%%%%%%%%%%%%%%%%%%%%%%%%%%%%%%%%%%%%%

\Fref{fig:bifInfCycle}(a1) shows a phase portrait of system~\eref{eq:sanDelay} showcasing the periodic orbit $\mathbf{\Gamma}$ (red curve) belonging to family $p_1$ at $\tau = -100$ and $a=-2.4$. The periodic orbit  $\mathbf{\Gamma}$ has period $T \approx 205.2721$ thereby satisfying $T/|\tau| \approx 2$. Notice in panel~(a1) that, unlike families of periodic orbits with scaling  $T/|\tau| \approx 1$, the periodic orbit $\mathbf{\Gamma}$ is not a TDS; that is, it does not approximate a homoclinic orbit. On the contrary, $\mathbf{\Gamma}$ is temporally localized in two points, indicated by asterisks in panel~(a1), in a way  that resembles a heteroclinic cycle. The corresponding temporal traces in each component are shown in panels~(a2) to (a4), where the two plateaus are clearly visible. However, after close inspection, one notices that these two points do not correspond to equilibria of system~\eref{eq:sanDelay}; therefore, the periodic orbit cannot be asymptotic to a heteroclinic cycle of the DDE.  Rather, the underlying limiting object is a period-two point of a reduced map, which is a well-known phenomenon in DDEs with large delay \cite{Mallet-Paret1986m}. To understand this concept, consider the time rescaling $s= \phi(t)= -t/\tau$ (for $\tau<0$) and $(X,Y,Z)(s)=(x,y,z)\circ\phi^{-1}(s)$ such that \eref{eq:sanDelay} takes the form
\begin{equation}\label{eq:sanDelayNewCoor} 
\begin{aligned} 
  \varepsilon X' &= aX + bY -aX^2+(\tilde{\mu}-\alpha Z)X(2-3X), \\[-0.5mm]
 \varepsilon Y' &= bX +aY -\frac32 bX^2-\frac32 aXY-2Y(\tilde{\mu}-\alpha
Z)+ \kappa \text{S}_{-1}\! \lp[ XY\rp], \\[-0.5mm]
 \varepsilon Z' &= cZ +\mu Z +\gamma XZ,
\end{aligned}
\end{equation}
where $\varepsilon=1/|\tau|$ is now a small parameter for large value of $|\tau|$. The left-hand side \eref{eq:sanDelayNewCoor} vanishes in the formal limit $\varepsilon=0$, which defines an  algebraic relation between $(X,Y,Z)(s)$ and $(X,Y,Z)(s+1)$. For  some cases, the derived algebraic relation defines a map, and, under certain conditions, period-two points of such a map have rigorously been shown  to generate families of periodic orbits with our desired localization and scaling properties for small $\varepsilon$ \cite{Mallet-Paret1986m}. This is not the case for  the formal limit of system~\eref{eq:sanDelayNewCoor}, since we can express $Y$ and $Z$ as
\begin{equation*}
\begin{aligned} 
  Y &= \dfrac{(-a c-2 c \tilde{\mu}) X +(a c-a \gamma -2 \alpha  \mu +3 c \tilde{\mu}-2 \gamma  \tilde{\mu}) X^2 + (a \gamma +3 \alpha  \mu +3 \gamma  \tilde{\mu}) X^3}{b (c + \gamma X)} \quad \text{and} \\[-0.5mm]
  Z &= \dfrac{\mu X}{c+ \gamma X},
\end{aligned}
\end{equation*}
thus obtaining a reduced algebraic relation between $X(s)$ and $X(s+1)$ given by 
\begin{equation} \label{eq:relationAtInfy}
    0 = bX +aY -\frac32 bX^2-\frac32 aXY-2Y(\tilde{\mu}-\alpha Z)+\kappa \text{S}_{-1}\lp[XY\rp],
\end{equation}
which can be interpreted as a  multivalued map in $X$. By expressing \eref{eq:relationAtInfy} with a common denominator, one obtains a zero problem in the form of a quartic polynomial on $X(s+1)$, which can be solved using symbolic algebra software. To find period-two points of Eq.~\eref{eq:relationAtInfy},
we solve the problem
\begin{equation*}
  \begin{aligned}
    0 &= bX +aY -\frac32 bX^2-\frac32 aXY-2Y(\tilde{\mu}-\alpha Z)+ \kappa \text{S}_{-1}\! \lp[XY\rp], \\[-0.5mm]
    0 &= \text{S}_{-1}\! \lp[ bX +aY -\frac32 bX^2-\frac32 aXY-2Y(\tilde{\mu}-\alpha Z) \rp]+ \kappa XY.
  \end{aligned}
\end{equation*}
For the parameter values chosen in panel~(a),  we find that the only two-period point with positive $X$-values is given by 
$$\lp(X(s) , X(s+1)\rp) \approx (0.6772,  0.275)$$
corresponding to $\mathbf{q_1}=(x,y,z) \approx (0.6772, 0.2083 ,0.2048)$ and $\mathbf{q_2} =(x,y,z) \approx (0.275, 0.195758, 0.0638)$ in the original coordinates. This is in good agreement with the values where the time traces are temporally localized in \fref{fig:bifInfCycle}(a2)-(a4). This result is indeed surprising as Eq.~\eref{eq:relationAtInfy} only defines a multivalued map.

%he formal limit of system~\eref{eq:sanDelayNewCoor} is not a map but the multivalued map:
%This is not the case for system~\eref{eq:sanDelayNewCoor} as the formal limit can be expressed as \ag{Insert equation here}

Now we focus on the limiting behavior of the homoclinic orbits along $\mathbf{H^-_t}$. We study this by starting from the high-periodic solution shown in panel~(a) and increasing $a$ towards $\mathbf{H^-_t}$. \Fref{fig:bifInfCycle}(b1) shows a phase portrait of an approximate homoclinic orbit. Notice that compared to panel~(a1), the homoclinic orbit is similar in shape and is  localized near the period-two points obtained by continuing $\mathbf{q_1}$ and $\mathbf{q_2}$ in $a$; however, the orbit has also started to be localized around $\mathbf{\widehat{0}}$ as is also confirmed by the temporal traces in panels~(b2)-(b4). Further numerical exploration (not shown) for large negative values of the delay past the values shown in \Fref{fig:bifInfCycle}(b) indicates that the temporal plateaus exhibited by the orbit do not exactly coincide with the values of $\mathbf{q_1}$ and $\mathbf{q_2}$ from the derived multivalued map~\eref{eq:relationAtInfy}, but actually have split up into smaller segments with nearby values. A careful inspection of these plateaus suggests the existence of a more complicated underlying object in  the multivalued map~\eref{eq:relationAtInfy}  to  which the homoclinic orbit along $\mathbf{H^-_t}$ is limiting as $\tau \rightarrow -\infty$. The exact nature of this transition is beyond the scope of this paper, and it will be an interesting avenue for future research as our numerics suggest that the generation of TDSs is also organized by periodic points and their bifurcations of reduced (multivalued) maps at infinity.

\section{Orbit flip case} \label{sec:OrbitFlipCase}
%%%%%%%%%%%%%%%%%%%%%%%%%%%%%%%%%%%%%%%%%%%%%%%%%%%%%%%%
\begin{table}[]
\centering
\begin{tabular}{|ccccccrcl|}
\hline
\multicolumn{5}{|c|}{\textbf{Parameters}} & \multicolumn{1}{c|}{\multirow{2}{*}{\textbf{\begin{tabular}[c]{@{}c@{}}Continuation\\ parameters\end{tabular}}}} & \multicolumn{3}{c|}{\multirow{2}{*}{\textbf{\begin{tabular}[c]{@{}c@{}}Stopping \\ condition\end{tabular}}}} \\ \cline{1-5}
\multicolumn{1}{|c|}{$a$} & \multicolumn{1}{c|}{$\mu$} & \multicolumn{1}{c|}{$\tilde{\mu}$} & \multicolumn{1}{c|}{$\kappa$} & \multicolumn{1}{c|}{$\tau$} & \multicolumn{1}{c|}{} & \multicolumn{3}{c|}{} \\ \hline
\multicolumn{9}{|c|}{\textbf{Homoclinic solution}} \\ \hline
\multicolumn{1}{|c|}{$-0.5$} & \multicolumn{1}{c|}{$0.0595$} & \multicolumn{1}{c|}{$0.03$} & \multicolumn{1}{c|}{$0$} & \multicolumn{1}{c|}{$0$} & \multicolumn{1}{c|}{$\kappa,\mu$} & $\mu$ & $=$ & $0.2$ \\ \hline
\multicolumn{1}{|c|}{$-0.5$} & \multicolumn{1}{c|}{$0.2$} & \multicolumn{1}{c|}{$0.03$} & \multicolumn{1}{c|}{$-0.4143$} & \multicolumn{1}{c|}{$0$} & \multicolumn{1}{c|}{$\kappa,\tau$} & $\kappa$ & $=$ & $-1$ \\ \hline
\multicolumn{1}{|c|}{$-0.5$} & \multicolumn{1}{c|}{$0.2$} & \multicolumn{1}{c|}{$0.03$} & \multicolumn{1}{c|}{$-1.0$} & \multicolumn{1}{c|}{$-0.5077$} & \multicolumn{1}{c|}{$\tilde{\mu},\tau$} & $\tilde{\mu}$ & $=$ & $0$ \\ \hline
\multicolumn{1}{|c|}{$-0.5$} & \multicolumn{1}{c|}{$0.2$} & \multicolumn{1}{c|}{$0$} & \multicolumn{1}{c|}{$-1.0$} & \multicolumn{1}{c|}{$-0.3927$} & \multicolumn{1}{c|}{\textemdash} & \multicolumn{3}{c|}{\textemdash} \\ \hline
\multicolumn{9}{|c|}{\textbf{Two-homoclinic solution}} \\ \hline
\multicolumn{1}{|c|}{$-0.7490$} & \multicolumn{1}{c|}{$0.0550$} & \multicolumn{1}{c|}{$0.03$} & \multicolumn{1}{c|}{$0$} & \multicolumn{1}{c|}{$0$} & \multicolumn{1}{c|}{$\kappa,\mu$} & $\mu$ & $=$ & $0.2$ \\ \hline
\multicolumn{1}{|c|}{$-0.7490$} & \multicolumn{1}{c|}{$0.2$} & \multicolumn{1}{c|}{$0.03$} & \multicolumn{1}{c|}{$-0.5613$} & \multicolumn{1}{c|}{$0$} & \multicolumn{1}{c|}{$\kappa,\tau$} & $\kappa$ & $=$ & $-1$ \\ \hline
\multicolumn{1}{|c|}{$-0.7490$} & \multicolumn{1}{c|}{$0.2$} & \multicolumn{1}{c|}{$0.03$} & \multicolumn{1}{c|}{$-1.0$} & \multicolumn{1}{c|}{$-0.4384$} & \multicolumn{1}{c|}{$a,\tau$} & $a$ & $=$ & $-0.5$ \\ \hline
\multicolumn{1}{|c|}{$-0.5$} & \multicolumn{1}{c|}{$0.2$} & \multicolumn{1}{c|}{$0.03$} & \multicolumn{1}{c|}{$-1.0$} & \multicolumn{1}{c|}{$-0.4122$} & \multicolumn{1}{c|}{$\tilde{\mu},\tau$} & $\tilde{\mu}$ & $=$ & $0$ \\ \hline
\multicolumn{1}{|c|}{$-0.5$} & \multicolumn{1}{c|}{$0.2$} & \multicolumn{1}{c|}{$0$} & \multicolumn{1}{c|}{$-1.0$} & \multicolumn{1}{c|}{$-0.2841$} & \multicolumn{1}{c|}{\textemdash} & \multicolumn{3}{c|}{\textemdash} \\ \hline
\multicolumn{9}{|c|}{\textbf{Period-doubling bifurcation}} \\ \hline
\multicolumn{1}{|c|}{$-0.6414$} & \multicolumn{1}{c|}{$0.0550$} & \multicolumn{1}{c|}{$0.03$} & \multicolumn{1}{c|}{$0$} & \multicolumn{1}{c|}{$0$} & \multicolumn{1}{c|}{$\kappa,\mu$} & $\mu$ & $=$ & $0.2$ \\ \hline
\multicolumn{1}{|c|}{$-0.6414$} & \multicolumn{1}{c|}{$0.2$} & \multicolumn{1}{c|}{$0.03$} & \multicolumn{1}{c|}{$-0.5613$} & \multicolumn{1}{c|}{$0$} & \multicolumn{1}{c|}{$\kappa,\tau$} & $\kappa$ & $=$ & $-1$ \\ \hline
\multicolumn{1}{|c|}{$-0.6414$} & \multicolumn{1}{c|}{$0.2$} & \multicolumn{1}{c|}{$0.03$} & \multicolumn{1}{c|}{$-1.0$} & \multicolumn{1}{c|}{$-0.4250$} & \multicolumn{1}{c|}{$a,\tau$} & $a$ & $=$ & $-0.5$ \\ \hline
\multicolumn{1}{|c|}{$-0.5$} & \multicolumn{1}{c|}{$0.2$} & \multicolumn{1}{c|}{$0.03$} & \multicolumn{1}{c|}{$-1.0$} & \multicolumn{1}{c|}{$7.8588$ (*)} & \multicolumn{1}{c|}{$\tilde{\mu},\tau$} & $\tilde{\mu}$ & $=$ & $0$ \\ \hline
\multicolumn{1}{|c|}{$-0.5$} & \multicolumn{1}{c|}{$0.2$} & \multicolumn{1}{c|}{$0$} & \multicolumn{1}{c|}{$-1.0$} & \multicolumn{1}{c|}{$7.8195$} & \multicolumn{1}{c|}{\textemdash} & \multicolumn{3}{c|}{\textemdash} \\ \hline
\multicolumn{9}{|c|}{\textbf{Saddle-node bifurcation of periodic orbits}} \\ \hline
\multicolumn{1}{|c|}{$-0.3852$} & \multicolumn{1}{c|}{$-0.0550$} & \multicolumn{1}{c|}{$-0.03$} & \multicolumn{1}{c|}{$0$} & \multicolumn{1}{c|}{$0$} & \multicolumn{1}{c|}{$\kappa, a$} & $a$ & $=$ & $-0.5$ \\ \hline
\multicolumn{1}{|c|}{$-0.5$} & \multicolumn{1}{c|}{$-0.0550$} & \multicolumn{1}{c|}{$-0.03$} & \multicolumn{1}{c|}{$-0.0050$} & \multicolumn{1}{c|}{$0$} & \multicolumn{1}{c|}{$\kappa,\tilde{\mu}$} & $\kappa$ & $=$ & $-1$ \\ \hline
\multicolumn{1}{|c|}{$-0.5$} & \multicolumn{1}{c|}{$-0.0550$} & \multicolumn{1}{c|}{$-0.1971$} & \multicolumn{1}{c|}{$-1.0$} & \multicolumn{1}{c|}{$0$} & \multicolumn{1}{c|}{$\tilde{\mu},\tau$} & $\tilde{\mu}$ & $=$ & $0$ \\ \hline
\multicolumn{1}{|c|}{$-0.5$} & \multicolumn{1}{c|}{$-0.0550$} & \multicolumn{1}{c|}{$0$} & \multicolumn{1}{c|}{$-1.0$} & \multicolumn{1}{c|}{$-0.9054$} & \multicolumn{1}{c|}{\textemdash} & \multicolumn{3}{c|}{\textemdash} \\ \hline
\end{tabular}\vspace{6pt}
\caption{Additional information for the computation of the bifurcation diagram for the orbit flip bifurcation of case~\textbf{B} in system~\eref{eq:sanDelay}.  Shown are the sequences of continuation runs to lift the homoclinic bifurcation and bifurcations of periodic solutions to non-zero $\tau$. Other parameters are given in \Sref{sec:NoDelaySan}.}\label{tab:orbitFlipCont}
\end{table}
%%%%%%%%%%%%%%%%%%%%%%%%%%%%%%%%%%%%%%%%%%%%%%%%%%%%%%%%

As mentioned in the Introduction and \cref{sec:NoDelaySan}, an orbit flip bifurcation is another mechanism that can create multi-pulse homoclinic and periodic solutions. In particular, we focus here on an orbit flip bifurcation of case~$\mathbf{B}$ which is the simplest case where a multi-pulse homoclinic solution arises; namely, a two-loop homoclinic solution. To compute the bifurcation diagram of this case in system~\cref{eq:sanDelay}, we proceed in the same fashion as in the resonant cases; more precisely, we first choose a homoclinic solution of system \cref{eq:san} at parameter values
$$(a,\mu,\tilde \mu)=(-0.5,0.0595,0.03),$$
and then lift the homoclinic solution (and periodic solutions of interest) to system~\cref{eq:sanDelay} in a series of two-parameter continuation runs which are described in \cref{tab:orbitFlipCont}. Notice in \cref{tab:orbitFlipCont} that we lift all solutions to $a=-0.5$ and $\tilde \mu =0$, such that the eigenvalues conditions of $\mathbf{\widehat{0}}$ satisfies the condition for case $\mathbf{B}$ \cite{Kisaka1993}. Unlike the resonant cases, we present the numerical bifurcation diagram of the orbit flip of case~$\mathbf{B}$ in the $(\tau,\mu)$-parameter plane; in this way, the eigenvalues of $\mathbf{\widehat{0}}$ do not change as parameters are varied and, when $\mu=0$, the homoclinic orbit will be contained in the $xy$-plane, see \cref{sec:intDelaySan}.

%%%%%%%%%%%%%%%%%%%%%%%%%%%%%%%%%%%%%%%%%%%%%%%%%%%%%%%%
\begin{figure}
\centering
\includegraphics{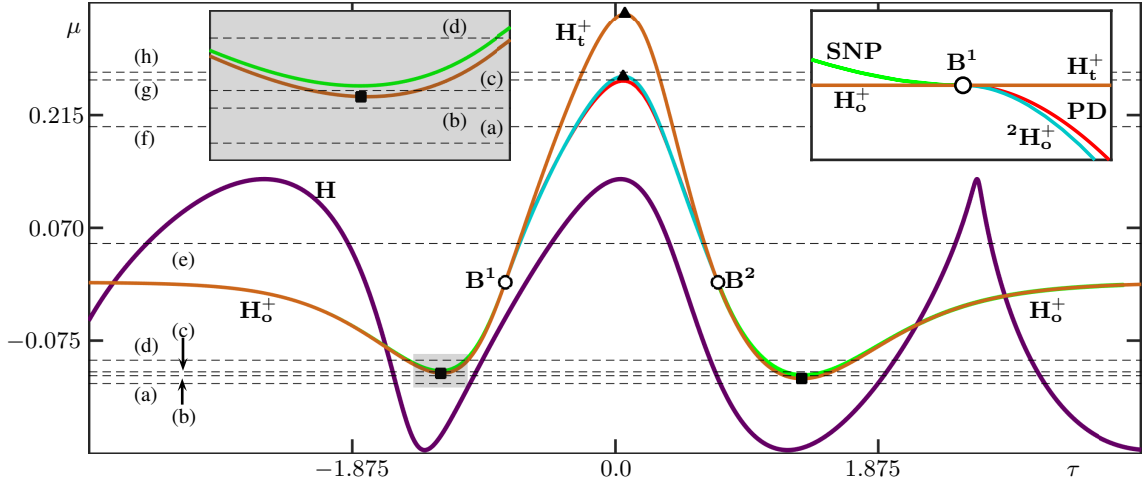}
\caption{Bifurcation diagram in the $(\tau, \mu)$-parameter plane of system~\eref{eq:sanDelay} near orbit flip bifurcation points $\mathbf{B^1}$ and $\mathbf{B^2}$ of case~\textbf{B}. The right inset shows an enlargement near the orbit flip bifurcation point $\mathbf{B^1}$ under a coordinate transformation to identify the homoclinic bifurcation curve with the horizontal axis. The left inset shows a enlargement of the grey region. Shown are the curve of homoclinic bifurcations~$\mathbf{H^+_o}$ and $\mathbf{H^+_t}$ (brown), saddle-node bifurcations $\mathbf{SNP}$ of periodic orbits (dark-green), period-doubling bifurcations $\mathbf{PD}$ (red),  two-homoclinic bifurcations $\mathbf{H^2_o}$ (cyan), and Hopf bifurcations $\mathbf{H}$ (purple curve).  The black square indicates the local minima $\mathbf{m[H_t^+]}$ in $\mu$ of $\mathbf{H_t^+}$, while the black triangles indicate maxima $\mathbf{M[H^+_t]}$  in $\mu$ of $\mathbf{H_t^+}$. The dashed curves indicate the slices considered in \fref{fig:bifDelayOrbitFlipSlice1} and \fref{fig:bifDelayOrbitFlipSlice2}. Here, $(a,\tilde{\mu},\kappa) \approx (-0.5,0.0,-1.0)$ and other parameters are given in \Sref{sec:NoDelaySan}.}\label{fig:bifDelayOrbitFlip}
\end{figure} 
%%%%%%%%%%%%%%%%%%%%%%%%%%%%%%%%%%%%%%%%%%%%%%%%%%%%%%%%

%%%%%%%%%%%%%%%%%%%%%%%%%%%%%%%%%%%%%%%%%%%%%%%%%%%%%%%%
\begin{figure}
\centering
\includegraphics{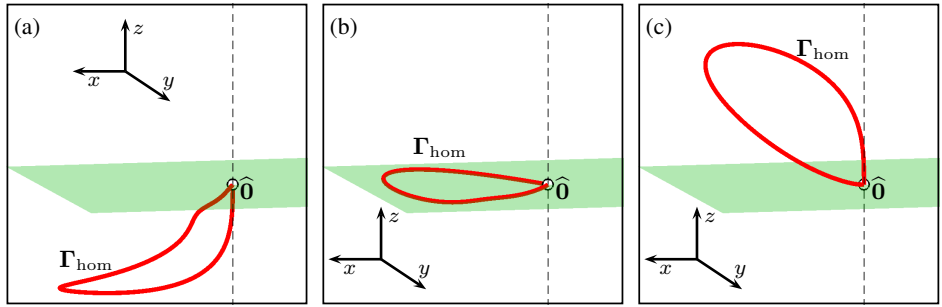}
\caption{Phase portraits of system~\eref{eq:sanDelay} showcasing the transition of the homoclinic orbit near the orbit flip bifurcation point $\mathbf{B^1}$ shown in \fref{fig:bifDelayOrbitFlip}.  Shown are the saddle equilibrium $\mathbf{\widehat{0}}$, the homoclinic orbit $\mathbf{\Gamma}_{\rm hom}$ (red curve), and the $xy$-plane (green plane). Panel~(a) corresponds to the situation along $\mathbf{H^+_{o}}$ for $(\mu, \tau) \approx  (-0.1000, -1.0492)$.  Panel~(b) corresponds to the situation at $\mathbf{B}^1$ for $(\mu, \tau) = (0.0000,-0.7840)$. Panel~(c) corresponds to the situation along $\mathbf{H^+_{t}}$ for $(\mu, \tau) =(0.1000, -0.5941)$.} \label{fig:orbitFlipTransition} 
\end{figure}
%%%%%%%%%%%%%%%%%%%%%%%%%%%%%%%%%%%%%%%%%%%%%%%%%%%%%%%%

\Fref{fig:bifDelayOrbitFlip} shows the bifurcation diagram in the $(\tau,\mu)$-parameter plane where there exists a single curve of homoclinic orbits (brown) which exists for all values of $\tau$, and exhibits two orbit flip bifurcations $\mathbf{B^1}$ and $\mathbf{B^2}$. By computing the period-doubling bifurcation curve $\mathbf{PD}$ (red) and two-homoclinic bifurcation $\mathbf{^2H^+_o}$ (cyan), we identify the segment of the homoclinic bifurcation curve joining $\mathbf{B^1}$ and $\mathbf{B^2}$ as non-orientable $\mathbf{H^+_t}$, and the other segments as orientable $\mathbf{H^+_o}$. The inset in \fref{fig:bifDelayOrbitFlip}  shows a magnification around the point $\mathbf{B^1}$ after a nonlinear transformation has been applied such that  $\mathbf{H^+_{o/t}}$ lies on the horizontal axis. In the inset, it is easily observed that the bifurcation curves emanating from $\mathbf{B^1}$ have the same topological organization as in the autonomous case presented in \fref{fig:bifAutoNoDelay}(a). Notice also the existence of a Hopf bifurcation curve $\mathbf{H}$ (purple) that exists for all $\tau$ values shown in the figure.

To illustrate and provide numerical evidence that the homoclinic orbits are indeed at an orbit flip configuration at $\mathbf{B^{1,2}}$, we show in \Fref{fig:orbitFlipTransition} numerically computed phase portraits of the homoclinic orbits (red)  along $\mathbf{H^+_{o/t}}$ before, at and after $\mathbf{B^{1}}$. In all panels, the $xy$-plane is represented by a green surface, and the $z$-axis is indicated as a dashed line. Panel~(a) shows the homoclinic orbit converging forward in time to $\mathbf{\widehat{0}}$ tangent to the $z$-axis (weak stable direction) below the $xy$-plane; here, the homoclinic orbit is orientable. On the other hand, the homoclinic orbit lies in the $xy$-plane at $\mathbf{B^{1}}$ as both $\mu$ and $\tilde \mu$ are zero, as shown in panel~(b). Here, the homoclinic orbit converges forward in time to $\mathbf{\widehat{0}}$ tangent to its strong stable direction which is contained in the $xy$-plane. By moving past this point, the homoclinic orbit becomes non-orientable and now converges forward in time to $\mathbf{\widehat{0}}$ tangent to the $z$-axis above the $xy$-plane.

%%%%%%%%%%%%%%%%%%%%%%%%%%%%%%%%%%%%%%%%%%%%%%%%%%%%%%%%
\begin{figure}
\centering
\includegraphics{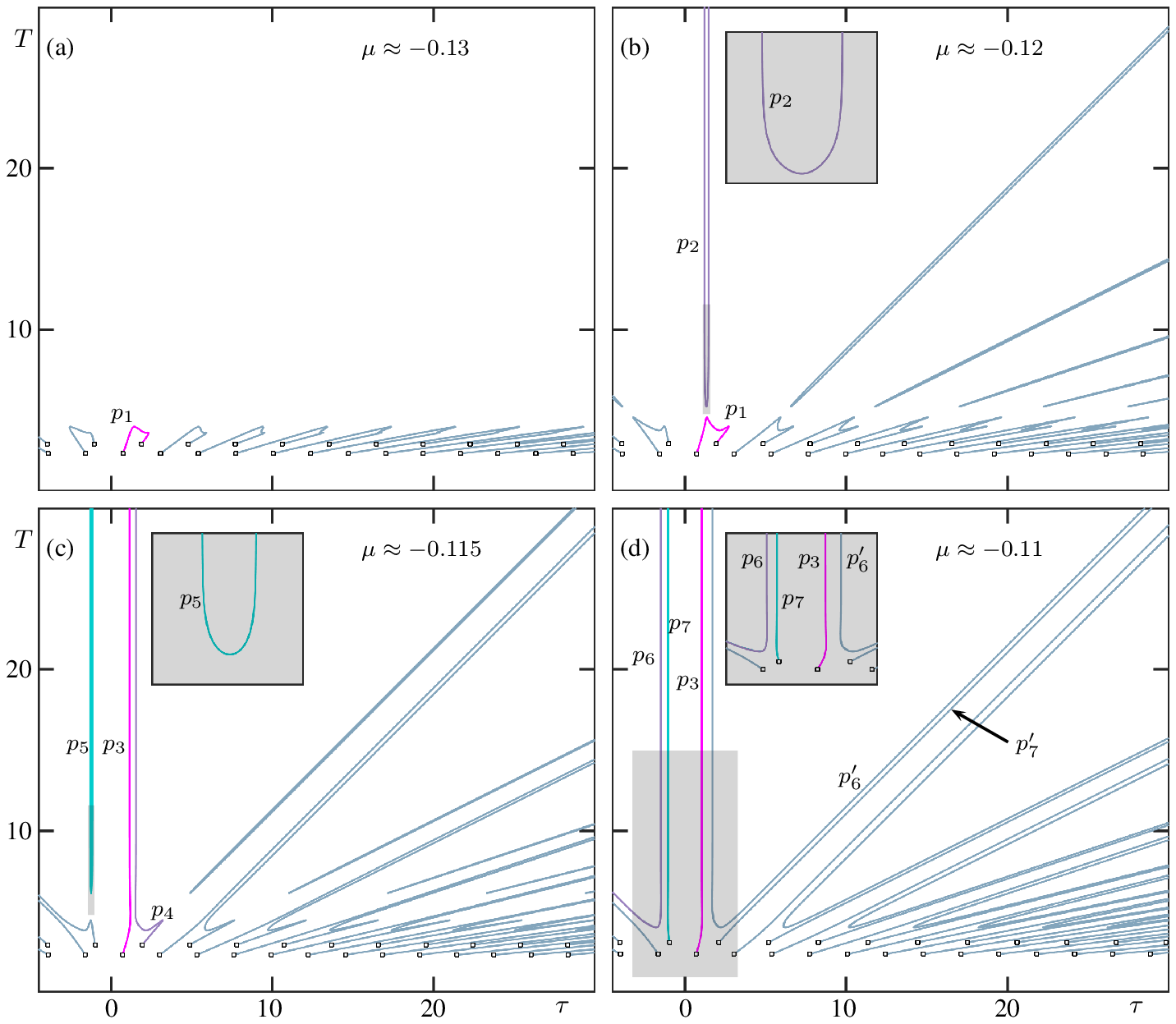}
\caption{One-parameter bifurcation diagrams in $\tau$ of system~\eref{eq:sanDelay} for different values of $\mu$ as indicated in \fref{fig:bifDelayOrbitFlip}. The value of $a$ for each slice is indicated on the top left side of each panel. Shown is the period $T$ of the periodic orbit with respect to $\tau$. The branch of periodic solutions emanating from the homoclinic bifurcation~$\mathbf{H^+_o}$ and $\mathbf{H^+_t}$ are colored cyan, magenta and light-purple, while the other periodic branches are colored gray.  The Hopf bifurcation points are indicated by small squares.} \label{fig:bifDelayOrbitFlipSlice1} 
\end{figure} 
%%%%%%%%%%%%%%%%%%%%%%%%%%%%%%%%%%%%%%%%%%%%%%%%%%%%%%%%

%%%%%%%%%%%%%%%%%%%%%%%%%%%%%%%%%%%%%%%%%%%%%%%%%%%%%%%%
\begin{figure}
\centering
\includegraphics{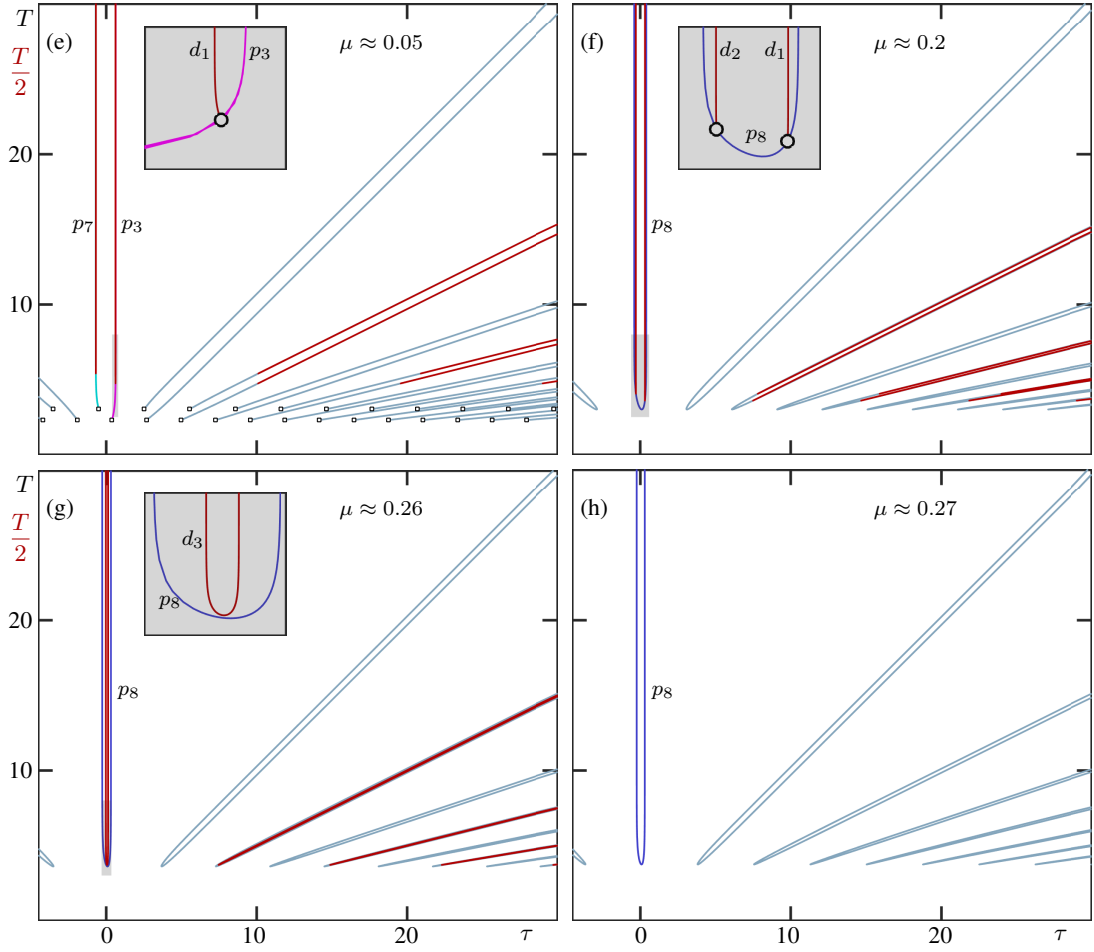}
\caption{Continued from \fref{fig:bifDelayOrbitFlipSlice1}. The branch of periodic solutions emanating the homoclinic bifurcation~$\mathbf{H_t}$ shown in  \fref{fig:bifDelayOriRes} is colored cyan, magenta and blue, and the ones emanating from the two-loop homoclinic bifurcation $\mathbf{^2H^+_o}$ and period-doubling bifurcation $\mathbf{PD}$ by red, while the other periodic branches are colored gray. The Hopf bifurcation points are indicated by small squares. } \label{fig:bifDelayOrbitFlipSlice2} 
\end{figure} 
%%%%%%%%%%%%%%%%%%%%%%%%%%%%%%%%%%%%%%%%%%%%%%%%%%%%%%%%

In the same spirit as the resonant cases, we now focus our attention on the existence of periodic orbits near and far away from the orbit flip points. The representative values of $\mu$  that were chosen are indicated by dashed lines in \fref{fig:bifDelayOrbitFlip}. \Fref{fig:bifDelayOrbitFlipSlice1}(a) shows the $\tau$-bifurcation diagram for $\mu \approx -0.13$, just below the occurrence of the homoclinic bifurcation curve $\mathbf{H^+_o}$. In this panel, a family of periodic orbits $p_1$ (magenta) exists and connects a pair of Hopf bifurcations points. Henceforth, the reappearing families share this property which corresponds to infinitely many intersections between the Hopf bifurcation curve $\mathbf{H}$ with this slice. As $\mu$ increases, we observe for positive $\tau$ a fold of homoclinic bifurcations $\mathbf{m[H^+_o]}$ in the extended system. Past this value, system~\eref{eq:sanDelay} exhibits two homoclinic bifurcation for fixed value of $\mu$ which is illustrated in \fref{fig:bifDelayOrbitFlipSlice1}(b) for $\mu\approx -0.12$. Here, a new family of periodic solutions $p_2$ (purple) is created, which has a global minimum in the period $T$ and connects both homoclinic bifurcations, as observed in panel~(b). As in the resonance cases, this global minimum maps close to saddle-node bifurcation that connects two branches of periodic solutions in the reappearing families. Increasing $\mu$ past the second fold of homoclinic bifurcation $\mathbf{m[H^+_o]}$ occurring for negative value $\tau$ and just below the saddle-node bifurcation $\mathbf{SNP}$, we encounter the situation showcased in panel~(c). Here, a new family of periodic solution $p_5$ (cyan), which connects the two new homoclinic bifurcations, emerges. Similarly to family $p_2$, the family $p_5$ also has a global period minimum, and its reappearing families have properties analogous to what was observed in $p_2$. Observe that, in-between panel~(b) and panel~(c), a fold of saddle-node bifurcations  $\mathbf{m[SNP]}$ in the extended system has occurred which allows $p_1$ and $p_2$ to interact; thus, creating the new branches $p_3$ (magenta) and $p_4$ (purple) in panel~(c). The same transition was also observed between \fref{fig:bifDelayOriResSlice1}(d) and \fref{fig:bifDelayOriResSlice1}(e) for the orientable resonant case.  The saddle-node bifurcation curve (green) associated with this transition is located in the direct vicinity of $\mathbf{H^+_o}$ in the $(\tau,\mu)$-parameter plane emanating from $\mathbf{B^2}$. This particular transition is also observed for family $p_5$ when $\mu$ is slightly increased to $-0.115$ past a second fold $\mathbf{m[SNP]}$ of saddle-node bifurcations in the extended system (occurring for negative values of $\tau$) along the saddle-node bifurcation curve $\mathbf{SNP}$ coming from $\mathbf{B^1}$ as shown in panel~(d). At the point $\mathbf{m[SNP]}$, the family $p_5$ interacts with one of the reappearing families of $p_4$ which affects the organization of these families. Indeed, after this transition, $p_5$ has split into the families $p_6$ (purple) and  $p_7$ (cyan); see panel~(d). Notice that as part of this reconfiguration, family $p_4$ has split into families $p'_6$ and $p'_7$ which correspond to the first reappearance of families $p_6$ and  $p_7$. Particularly, the family $p_7$  connects the homoclinic and Hopf bifurcation points that occur for small negative values of $\tau$. On the other hand, the family $p_6$ arises from a homoclinic bifurcation point; and, following the family to large negative of $\tau$, one notices that its period grows beyond bound scaling linearly with $|\tau|$. This feature is quite special from what we have seen so far, as the family $p'_6$ (first reappearance of $p_6$) also possesses this feature. 

Now we focus our attention on the situation past the bifurcation point $\mathbf{B^{1,2}}$; thus, we show the $\tau$-bifurcation diagram for  $\mu \approx 0.05$ in \fref{fig:bifDelayOrbitFlipSlice2}(e). We remark that the slice shown in panel~(e) is above the horizontal asymptotes of $\mathbf{H^+_o}$ in $\tau$, see \fref{fig:bifDelayOrbitFlip}. As consequence, families $p_6$ and $p'_6$ have disappeared in panel~(e). This transition is possibly related to an object in the multivalued map~\eref{eq:relationAtInfy}, similar to the non-orientable resonant homoclinic case studied in \Sref{sec:twoInfty}. We remark that families $p_6$ and $p'_6$ are related by the reappearance map; thus, both must coexist at a fixed $a$-slice. For this reason, it is clear that the horizontal asymptotes of $\mathbf{H^+_o}$ for positive and negative $\tau$ coincide. This transition is of interest to the organization of TDS in system~\eref{eq:sanDelay}; however, its exact nature regarding the multivalued map~\eref{eq:relationAtInfy} is beyond the scope of this paper, and it is an interesting avenue for future research.

We now concentrate on the families $d_1$ (red) and $d_2$ (red) of period-doubled orbits that are shown in panel~(e). Recall that system~\eref{eq:sanDelay} exhibits two period-doubling bifurcations and two two-homoclinic bifurcations at $\mu \approx 0.05$. This manifest in the $\tau$-bifurcation diagram in the following way: each of these two families connects a two-homoclinic bifurcation with a period-doubling bifurcation along $p_3$ or $p_7$. Notice that as we increase $\mu$, the extended system~\eref{eq:sanDelayExt} exhibits a fold of Hopf bifurcations $\mathbf{M[H]}$ corresponding to the maximum $\mu$-value along the Hopf bifurcation curve $\mathbf{H}$. At $\mathbf{M[H]}$, pairs of Hopf bifurcation observed in panel~(e) coincide and vanish as $\mu$ increases. Indeed, when we increase $\mu$ past this point, one sees in panel~(f) that the two periodic families $p_3$ and $p_7$ connect up to create the family $p_8$ (blue) which now has a global period minimum. As a result, family $p_8$ now exhibits both period-doubling bifurcations from which $d_1$ and $d_2$ emanate without creating topological changes to these branches. Notice that, by increasing $\mu$ further, we observe a fold $\mathbf{M[PD]}$ of period-doubling bifurcations in the extended system~\eref{eq:sanDelayExt} corresponding to the maximum $\mu$-value along the period-doubling bifurcation curve in the $(\tau,\mu)$-parameter plane.   Here, the two period-doubling bifurcation points coincide and disappear as $\mu$ is increased. Panel~(g) shows the situation at $\mu \approx 0.26$ after this transition.   Here, the period-doubled families $d_1$ and $d_2$ have merged, and we refer to the resulting family as $d_3$. The family $d_3$ now has a global period  minimum, and it is not connected to family $p_8$ through any bifurcation. Notice that the families that reappear from $d_3$ now exhibit saddle-node bifurcations nearby the $\tau$ values where the minimum along $d_3$ is mapped.  We remark that in between panel~(f) and (g), either the period-doubling bifurcations along the family $p_8$ must have changed criticality, or the families $d_1$ and $d_2$ must exhibit saddle-node bifurcations which undergo a fold of saddle-node bifurcations in the extended system (saddle transition). The exact nature of this particular transition involves fold bifurcations of period-doubled orbits which are beyond the scope of this paper. Finally, the period-minimum of the family $p_3$ increases without bounds as $\mu$ nears the fold of two-homoclinic bifurcations $\mathbf{M[^2H^+_o]}$  in the extended system, which corresponds to the maximum $\mu$-value along the curve $\mathbf{^2H^+_o}$; see \fref{fig:bifDelayOrbitFlip}. Panel~(h) shows the situation for $\mu \approx 0.27$ past the point $\mathbf{M[^2H^+_o]}$, where the family $p_3$ of period-doubled orbits has vanished from the $\tau$-bifurcation diagram, marking the disappearance of the bound two-pulse TDS.

\section{Conclusions} \label{sec:conclusions}
We have identified the link between codimension-two homoclinic bifurcations to a real saddle and the existence of bound two-pulse temporal dissipative solitons (TDSs) in delay differential equations (DDEs). Namely, we studied the orientable and non-orientable case of resonant bifurcations \cite{Chow1990}, and the orbit flip bifurcation of case~\textbf{B} \cite{Sandstede1993a}. To accomplish this, we have extended Sandstede's model~\cite{san1}  with a time-shift parameter, and  have \textemdash for the first time\textemdash showcased the numerical unfolding of the aforementioned codimension-two homoclinic bifurcations in the presence of time-delay. This was achieved by employing a numerical homotopy step from solutions with zero time-delay, obtained with \textsc{Auto07p}, and lifting them to non-zero delay with DDE-BIFTOOL. As a secondary result, we have also identified the role of minima and maxima along homoclinic bifurcation curves in the two-parameter plane, as mechanisms to create and destroy families of TDSs. Furthermore, we have reported on signatures of the dynamics of a multivalued-map derived in the limit $\tau \rightarrow -\infty$ in the organization of TDSs.

This work has primarily focussed on the existence of TDS, leaving the stability analysis as an interesting avenue for future research \cite{Yanchuk2019}. Homoclinic bifurcations have been studied extensively in ODEs \cite{Homburg2010, Leonid1}, and future work will be dedicated to studying the relationship between TDSs and homoclinic bifurcations further. Candidates include orbit flip/inclination flip bifurcations of case~\textbf{C}, which give rise to  $n$-homoclinic orbits \cite{Hom1,Hupkes2009, Kisaka1993}. It has been shown that such bifurcations can also create exotic pulsing behavior organized by bifurcations with objects at infinity \cite{Giraldo2020i}. Another interesting avenue will be TDSs with oscillating tails, which can be shown to correspond to Shilnikov bifurcations by the traveling wave ansatz for TDS \cite{Yanchuk2019}. We remark that homoclinic bifurcations are not the only mechanism that can organize TDSs, as heteroclinic cycles have been shown to generate interesting families of TDSs with particular signatures and solutions with multiple plateaus \cite{Ruschel2020, Ruschel2019d}.
Even more generally, the traveling wave equation for TDS does not only pertain DDEs but also traveling pulses (and also waves) in coarse-grained spatially-extended or coupled systems \cite{Hupkes2020}, placing our work as part of a comprehensive program to understand the generation of TDSs in applications.

\section*{Acknowledgments}
Andrus Giraldo was supported by KIAS Individual Grant No. CG086101 at Korea Institute for Advanced Study. We thank Bernd Krauskopf, Matthias Wolfrum and Serhiy Yanchuk for helpful comments and insightful discussions. 

\bibliographystyle{siam} 
% using 2 bib files for now % Stefan 08.05.22
\bibliography{GR_ResFlipDelay_arXiv}

\appendix

\section{Proof of \Cref{thm:reap-bif}}\label{sec:appendix}
Consider Eq.~\eref{eq:DDE} with sufficiently smooth right hand side $f: \C^d\times \C^d \to \C^d, (x,y)\mapsto f(x,y),$ for $d\in\mathbb{N}$.  Let $\tilde u$ be a periodic solution of Eq.~\eref{eq:DDE}, at $\tau = \tilde{\tau}$ and period $T$, with  resonant Floquet multiplier $\mu$ as in \Cref{thm:reap-bif} and respective eigenfunction $w\in C(\R,\mathbb{C}^d)$ satisfying
\begin{equation}
	\dot{w}(t)=A(t)w(t) + B(t)w(t- \tilde\tau),\quad  \mu w(t)=w(t+T), \label{eq:def-w-tildetau}
\end{equation} 
where $A,B:\R\to\C^{d\times d}$ are the $T$-periodic matrix coefficient functions $$ A(t)=\partial_{x}f(\tilde{u}(t),\tilde{u}(t- \tilde\tau)) \mbox{ and } B(t) = \partial_{y}f(\tilde{u}(t),\tilde{u}(t- \tilde\tau)).$$

Let $k$ be an integer and $\tau_k := \tilde \tau +jkT$, where $j$ is such that $\mu^j=1$. By using the reappearance map~\eref{eq:reap}, it follows that $\tilde u$ is a periodic solution of Eq.~\eref{eq:DDE} for $\tau_k$. It follows from the Floquet condition of \eref{eq:def-w-tildetau} that $$w(t- \tilde \tau)= w(t- \tau_k + jkT) = \mu^{jk}w(t- \tau_k)=w(t- \tau_k),$$ and Eq.~\eref{eq:def-w-tildetau} can be recast as
\begin{equation}\label{eq: def-w-tauk}
  \dot{w}(t)  =A(t)w(t)+ B(t)w(t- \tau_k)  ,\quad  \mu w(t)=w(t+T).
\end{equation}
Thus, $(\mu,w)$ is a solution of the Floquet problem \eref{eq: def-w-tauk} for $\tilde{u}$ at $\tau=\tau_k$.

\end{document}